\documentclass[a4paper,12pt,amsart,frenchb]{article}

\usepackage{amsmath,amsbsy,amsfonts,amssymb}
\usepackage[french]{babel}
\newcommand{\ass}[2]{\vskip0.3cm\noindent
{\bf {#1}}. { \sl {#2}}\vskip0.3cm\noindent
}
  
\begin{document}

 \title{  Les facteurs de transfert pour les groupes classiques: un formulaire}
\author{J.-L. Waldspurger}
\date{3 novembre 2009}
\maketitle

{\bf Introduction}

Le facteur de transfert est un terme qui intervient de fa\c{c}on cruciale dans la  stabilisation de la formule des traces. Il a \'et\'e d\'efini par Langlands et Shelstad dans le cas de l'endoscopie ordinaire ([LS]), puis dans le cadre plus g\'en\'eral de l'endoscopie tordue par Kottwitz et Shelstad ([KS]).  La d\'efinition est un peu abstraite et c'est un exercice amusant de l'expliciter dans des cas particuliers.  On calcule ici les facteurs de transfert pour les groupes classiques. Pr\'ecis\'ement pour les groupes symplectiques, sp\'eciaux orthogonaux, unitaires, ainsi que pour les groupes lin\'eaires tordus et pour les groupes tordus d\'eduits  de groupes unitaires par changement de base. Le point  de vue adopt\'e, comme dans [W] chapitre X, est de param\'etrer les classes de conjugaison, ou de conjugaison stable, des \'el\'ements semi-simples r\'eguliers d'un tel groupe par des donn\'ees de g\'eom\'etrie \'el\'ementaire. Par exemple, la classe de conjugaison stable  d'un tel \'el\'ement est d\'etermin\'ee par les
valeurs propres de l'\'el\'ement dans la repr\'esentation naturelle du groupe. On peut expliciter les facteurs de transfert \`a l'aide de ces param\`etres. Le r\'esultat est la proposition 1.10 ci-dessous. Signalons que Kottwitz a trouv\'e, au moins dans le cas des alg\`ebres de Lie, une autre fa\c{c}on de calculer ces facteurs, \`a l'aide de sections de Kostant. Cette m\'ethode s'est  av\'er\'ee particuli\`erement fructueuse. On esp\`ere n\'eanmoins que les formules pr\'esent\'ees ici pourront trouver quelques applications. 
Dans la premi\`ere section, on d\'efinit  les groupes consid\'er\'es, les param\`etres utilis\'es, les $L$-groupes et les donn\'ees endoscopiques. A propos de ces derni\`eres, signalons qu'il importe de les fixer pr\'ecis\'ement. En effet, le facteur de transfert d\'epend  vraiment des donn\'ees endoscopiques, et pas seulement de leur classe d'\'equivalence (on s'en convainc en consid\'erant les donn\'ees endoscopiques de $GL_{1}$). Ce point n'est peut-\^etre pas assez soulign\'e dans la litt\'erature. La section se termine par l'\'enonc\'e du r\'esultat. Les d\'emonstrations n'ont aucun int\'er\^et: il s'agit simplement d'expliciter les d\'efinitions dans chacun des cas. Toutefois, il para\^{\i}trait peut-\^etre bizarre de ne pr\'esenter aucune preuve. Dans la deuxi\`eme section, on a donc r\'edig\'e une preuve: celle du cas des groupes lin\'eaires tordus, en dimension impaire.

\bigskip
 
 \section{ Les r\'esultats}
 
 \bigskip
 
 \subsection{Notations}
 
 Soit $F$ un corps local de caract\'eristique nulle. On suppose que $F$ n'est pas le corps des complexes, la th\'eorie \'etant triviale dans ce cas. On fixe une cl\^oture alg\'ebrique $\bar{F}$ de $F$. On appelle ici extension alg\'ebrique de $F$ un sous-corps de $\bar{F}$ contenant $F$.  On note $Gal(\bar{F}/F)$, resp. $W_{F}$, le groupe de Galois, resp. le groupe de Weil, de $\bar{F}/F$. Par l'homomorphisme $W_{F}\to F^{\times}$ du corps de classes, tout homomorphisme $\chi$ de $F^{\times}$ d\'efinit un homomorphisme de $W_{F}$ que l'on note encore $\chi$. Si $E$ est une extension quadratique de $F$, on note $\tau_{E/F}$ l'unique \'el\'ement non trivial du groupe de Galois $Gal(E/F)$ et $sgn_{E/F}$ le caract\`ere quadratique de $F^{\times}$ dont le noyau est le groupe des normes $Norm_{E/F}(E^{\times})$.
 
 Si $G$ est un groupe r\'eductif d\'efini sur $F$, on note $\mathfrak{g}$ son alg\`ebre de Lie.
 \bigskip
 
 \subsection{D\'efinition des groupes et groupes tordus}
 
 On \'etudiera dans la suite diff\'erents cas pour lesquels on introduit ici les notations de base, que l'on ne rappellera plus.
 
 {\bf Le cas symplectique.} Soit $V$ un espace vectoriel sur $F$ de dimension finie $d$, muni d'une forme symplectique $q$ ($d$ est donc pair). On note $G$ le groupe symplectique de $(V,q)$.
 
{\bf Le cas sp\'ecial orthogonal.} Soit $V$ un espace vectoriel sur $F$ de dimension finie $d$, muni d'une forme quadratique (c'est-\`a-dire une forme bilin\'eaire sym\'etrique) non d\'eg\'en\'er\'ee $q$. On note $G$ le groupe sp\'ecial orthogonal de $(V,q)$. Ce cas, comme certains des cas suivants, se subdivise en un cas pair et un cas impair selon la parit\'e de $d$. Dans le cas pair, on d\'efinit le discriminant  $\delta$ de $q$ par $\delta=(-1)^{d/2}det(q)$ (par convention, $\delta=1$ si $d=0$). C'est un \'el\'ement de $F^{\times}/F^{\times,2}$. On l'a normalis\'e de sorte que $\delta=1$ si $(V,q)$ est somme orthogonale de plans hyperboliques. On exclut le cas o\`u $d=2$ et $\delta=1$.

 {\bf Le cas du groupe lin\'eaire tordu.} Soit $V$ un espace vectoriel sur $F$ de dimension finie $d$. On note $G$ le groupe (alg\'ebrique)des automorphismes lin\'eaires de $V$. On note $\tilde{G}$ l'ensemble (ou plus exactement la vari\'et\'e alg\'ebrique) des formes bilin\'eaires non d\'eg\'en\'er\'ees sur $V\times V$. Le groupe $G$ agit \`a gauche et \`a droite sur $\tilde{G}$ de la fa\c{c}on suivante. Soient $g,g'\in G$ et $\tilde{x}\in \tilde{G}$. Alors $g\tilde{x}g'$ est la forme bilin\'eaire
 $$(v,v')\mapsto \tilde{x}(g^{-1}v,g'v').$$
 Ces actions font de $(G,\tilde{G})$ un "groupe tordu", dans la terminologie de Labesse. 
 
 Consid\'erons le groupe lin\'eaire $GL_{d}$. On note $\theta_{d}$ son automorphisme d\'efini par
 $\theta_{d}(g)=J_{d}{^tg}^{-1}J_{d}^{-1}$, o\`u $J_{d}$ est la matrice antidiagonale dont les coefficients non nuls sont d\'efinis par $(J_{d})_{k,d+1-k}=(-1)^{k}$. Le carr\'e de $\theta_{d}$ est l'identit\'e et on peut introduire le produit semi-direct $GL_{d}^+=GL_{d}\rtimes \{1,\theta_{d}\}$. Fixons une base $(e_{k})_{k=1,...,d}$ de $V$. Elle nous permet d'identifier $G$ \`a $GL_{d}$. Fixons de plus un \'el\'ement $\nu\in F^{\times}$. Notons $\tilde{\theta}$ l'\'el\'ement de $\tilde{G}(F)$ d\'efini par les \'egalit\'es
 $$\tilde{\theta}(e_{k},e_{l})=\nu(-1)^{k}\delta_{k,d+1-l},$$
 o\`u le dernier terme est le symbole de Kronecker. On peut identifier $\tilde{G}$ \`a la composante non neutre $GL_{d}\theta_{d}$ de $GL_{d}^+$ par l'application $g\tilde{\theta}\mapsto g\theta_{d}$ pour tout $g\in G$.
 
 {\bf Le cas unitaire.}  Soit $E$ une extension quadratique  de $F$. Soit $V$ un espace vectoriel sur $E$ de dimension $d$, muni d'une forme hermitienne non d\'eg\'en\'er\'ee $q$. Pr\'ecisons notre convention concernant la sesquilin\'earit\'e: on a $q(\lambda v,\lambda' v')=\tau_{E/F}(\lambda)\lambda' q(v,v')$ pour tous $v,v'\in V$, $\lambda,\lambda'\in E$.  On note $G$ le groupe unitaire de $(V,q)$.
 
 {\bf Le cas du changement de base du groupe unitaire.} Soient $E$ une extension quadratique  de $F$ et $V$ un espace vectoriel sur $E$ de dimension $d$. On note $G$ le groupe alg\'ebrique sur $F$, obtenu par restriction des scalaires, tel que $G(F)$ soit le groupe des automorphismes $E$-lin\'eaires de $V$. On note $\tilde{G}$ l'ensemble (ou plus exactement la vari\'et\'e alg\'ebrique sur $F$) des formes sesquilin\'eaires non d\'eg\'en\'er\'ees sur $V$. Le groupe $G$ agit \`a gauche et \`a droite sur $\tilde{G}$ par la m\^eme formule que dans le cas du groupe lin\'eaire tordu. Ces actions font de $(G,\tilde{G})$ un "groupe tordu". 
 
 Notons $\theta_{d,E/F}$ l'automorphisme de $GL_{d}(E)$ d\'efini par $\theta_{d,E/F}(g)=J_{d}\tau_{E/F}(^tg^{-1})J_{d}^{-1}$, o\`u $\tau_{E/F}(^tg^{-1})$ est obtenu en appliquant $\tau_{E/F}$ aux coefficients de $^tg^{-1}$. Introduisons le produit semi-direct $GL_{d,E/F}^+(F)=GL_{d}(E)\rtimes \{1,\theta_{d,E/F}\}$. Fixons une base $(e_{k})_{i=1,...,d}$ de $V$ sur $E$. Elle nous permet d'identifier $G(F)$ \`a $GL_{d}(E)$. Fixons de plus un \'el\'ement $\nu\in E^{\times}$. Notons $\tilde{\theta}$ l'\'el\'ement de $\tilde{G}(F)$ d\'efini par les \'egalit\'es
 $$\tilde{\theta}(e_{k},e_{l})=\nu(-1)^{k}\delta_{k,d+1-l}.$$
 On peut identifier $\tilde{G}(F)$ au sous-ensemble $GL_{d}(E)\theta_{d,E/F}$ de $GL_{d,E/F}^+(F)$ par l'application $g\tilde{\theta}\mapsto g\theta_{d,E/F}$ pour tout $g\in G(F)$. On a d\'ecrit ainsi une identification des ensembles de points \`a valeurs dans $F$ mais il est facile de l'alg\'ebriser.
 
 \bigskip
 
 \subsection{Classes de conjugaison d'\'el\'ements semi-simples  suffisamment r\'eguliers}
 
 Hormis le cas sp\'ecial orthogonal pair,  "suffisamment" r\'egulier signifie pour nous fortement r\'egulier. Dans le cas sp\'ecial orthogonal pair, les \'el\'ements qui ont des valeurs propres $\pm 1$ peuvent \^etre fortement r\'eguliers, tout en n'ayant pas les m\^emes propri\'et\'es que les \'el\'ements en position g\'en\'erale. On entend alors par \'el\'ement semi-simple  suffisamment r\'egulier un \'el\'ement  semi-simple fortement r\'egulier qui n'a aucune valeur propre $\pm 1$.

 {\bf Le cas symplectique.} Donnons-nous la collection d'objets suivante:
 
 $\bullet$ un ensemble fini $I$;
 
 $\bullet$ pour tout $i\in I$, une extension finie $F_{\pm i}$ de $F$ et une $F_{\pm i}$-alg\`ebre commutative $F_{i}$ de dimension $2$ sur $F_{\pm i}$ (c'est-\`a-dire, ou bien $F_{i}$ est une extension quadratique de $F_{\pm i}$, ou bien $F_{i}=F_{\pm i}\oplus F_{\pm i}$); on note $\tau_{i}$ l'unique automorphisme non trivial de $F_{i}/F_{\pm i}$;
 
 $\bullet$ pour tout $i\in I$, un \'el\'ement $c_{i}\in F_{i}^{\times}$ tel que $\tau_{i}(c_{i})=-c_{i}$ et un \'el\'ement $x_{i}\in F_{i}^{\times}$ tel que $x_{i}\tau_{i}(x_{i})=1$.
 
 On suppose que $d=\sum_{i\in I}[F_{i}:F]$. Posons $W=\oplus_{i\in I}F_{i}$ et d\'efinissons une forme symplectique $q_{W}$ sur $W$ par 
 $$(1) \qquad q_{W}(\sum_{i\in I}w_{i},\sum_{i\in I}w'_{i})=\sum_{i\in I}trace_{F_{i}/F}(\tau_{i}(w_{i})w'_{i}c_{i}).$$
  
 Fixons un isomorphisme de $(W,q_{W})$ sur $(V,q)$. Notons $x$ l'\'el\'ement de $G(F)$ qui, modulo cet isomorphisme, est d\'efini par l'\'egalit\'e:
 $$(2) \qquad x(\sum_{i\in I} w_{i})=\sum_{i\in I}x_{i}w_{i}.$$
  Cet \'el\'ement est semi-simple et sa classe de conjugaison par $G(F)$ ne d\'epend pas de l'isomorphisme choisi. On d\'ecrit facilement \`a quelles conditions cet \'el\'ement est  r\'egulier. Disons simplement que si la famille $(x_{i})_{i\in I}$ est "en position g\'en\'erale", $x$ est   r\'egulier (ce qui \'equivaut \`a fortement r\'egulier dans le cas symplectique). Inversement toute classe de conjugaison semi-simple  r\'eguli\`ere dans $G(F)$ est obtenue par ce proc\'ed\'e. Donc, \`a tout \'el\'ement semi-simple  r\'egulier $x\in G(F)$, on peut associer des donn\'ees comme ci-dessus. L'ensemble $I$, les extensions $F_{\pm i}$ et $F_{i}$ et les \'el\'ements $x_{i}$ sont presque uniquement d\'etermin\'es (on peut \'evidemment remplacer $I$ par un autre ensemble de m\^eme nombre d'\'el\'ements, et chaque triplet $(F_{\pm i},F_{i},x_{i})$ par un triplet qui lui est isomorphe sur $F$). Par contre, les $c_{i}$ ne sont d\'etermin\'es qu'\`a multiplication pr\`es par le groupe des normes $Norm_{F_{i}/F_{\pm i}}(F_{i}^{\times})$.
  
  \bigskip
 
 {\bf Le cas sp\'ecial orthogonal impair}. On consid\`ere une collection d'objets comme dans le cas pr\'ec\'edent, \`a ceci pr\`es que l'on suppose maintenant que $\tau_{i}(c_{i})=c_{i}$. On suppose que $d=1+\sum_{i\in I}[F_{i}:F]$. On construit encore l'espace $W$, la forme $q_{W}$ qui est maintenant  quadratique et on suppose qu'il existe un espace $D$ de dimension $1$ muni d'une forme bilin\'eaire sym\'etrique et non d\'eg\'en\'er\'ee $q_{D}$ de sorte que $(W\oplus D,q_{W}\oplus q_{D})$ soit isomorphe \`a $(V,q)$. On fixe un tel isomorphisme. On introduit l'\'el\'ement $x\in G(F)$ qui, modulo cet isomorphisme, agit par la formule (2) sur $W$ et par l'identit\'e sur $D$. On a alors les m\^emes propri\'et\'es que dans le cas symplectique.
 
 \bigskip
 
 {\bf Le cas sp\'ecial orthogonal pair.} On consid\`ere une collection d'objets comme dans le cas symplectique, \`a ceci pr\`es que l'on suppose, comme dans le cas sp\'ecial orthogonal impair, que $\tau(c_{i})=c_{i}$. On suppose que $d=\sum_{i\in I}[F_{i}:F]$. On construit l'espace $W$ et la forme quadratique $q_{W}$. On suppose $(W,q_{W})$ isomorphe \`a $(V,q)$ et on fixe un isomorphisme. On construit  l'\'el\'ement $x\in G(F)$ comme dans le cas symplectique. On a essentiellement les m\^emes propri\'et\'es que dans le cas symplectique. Il y a toutefois un changement. C'est seulement la classe de conjugaison de $x$ par le groupe orthogonal de $V$ qui est bien d\'etermin\'ee. Or, pour un \'el\'ement semi-simple suffisamment r\'egulier de $G(F)$, sa classe de conjugaison par le groupe orthogonal se d\'ecompose en deux classes de conjugaison par $G(F)$. Autrement dit, notre collection d'objets $I$, $(F_{i})_{i\in I}$ etc... param\`etre non pas des classes de conjugaison par $G(F)$ mais des couples de telles classes.
 
 \bigskip
 
 {\bf Le cas du groupe lin\'eaire tordu, avec $d$ pair.} Donnons-nous la collection d'objets suivante:
 
 $\bullet$ un ensemble fini $I$;
 
 $\bullet$ pour tout $i\in I$, une extension finie $F_{\pm i}$ de $F$ et une $F_{\pm i}$-alg\`ebre commutative $F_{i}$ de dimension $2$ sur $F_{\pm i}$;
 
 $\bullet$ pour tout $i\in I$, un \'el\'ement $x_{i}\in F_{i}^{\times}$.
 
 On suppose $d=\sum_{i\in I}[F_{i}:F]$. On fixe un isomorphisme de $V$ sur $\oplus_{i\in I}F_{i}$. Modulo cet isomorphisme, on d\'efinit un \'el\'ement $\tilde{x}\in \tilde{G}(F)$ par la formule
 $$\tilde{x}(\sum_{i\in I}w_{i},\sum_{i\in I}w'_{i})=\sum_{i\in I}trace_{F_{i}/F}(\tau_{i}(w_{i})w'_{i}x_{i}).$$
 La classe de conjugaison par $G(F)$ de cet \'el\'ement $\tilde{x}$ est bien d\'etermin\'ee. Si la famille $(x_{i})_{i\in I}$ est "en position g\'en\'erale", $\tilde{x}$ est fortement r\'egulier. Inversement toute classe de conjugaison semi-simple fortement r\'eguli\`ere dans $\tilde{G}(F)$ est obtenue par ce proc\'ed\'e. Autrement dit, on peut associer \`a une telle classe des param\`etres $I$, $(F_{i})_{i\in I}$ etc... L'ensemble $I$ et les extensions $F_{\pm i}$ et $F_{i}$ sont d\'etermin\'es de fa\c{c}on essentiellement unique. Par contre les $x_{i}$ ne le sont qu'\`a multiplication pr\`es par le groupe $Norm_{F_{i}/F_{\pm i}}(F_{i}^{\times})$. 
 
 \bigskip
 
{\bf Le cas du groupe lin\'eaire tordu, avec $d$ impair.} On se donne une collection d'objets comme dans le cas pr\'ec\'edent, plus un \'el\'ement $x_{D}\in F^{\times}$. On pose $D=F$ et on fixe un isomorphisme de $V$ sur $D\oplus (\oplus_{i\in I}F_{i})$. Modulo cet isomorphisme, on d\'efinit un \'el\'ement $\tilde{x}\in \tilde{G}(F)$ par la formule
 $$\tilde{x}(w_{D}+\sum_{i\in I}w_{i},w'_{D}+\sum_{i\in I}w'_{i})=x_{D}w_{D}w'_{D}+\sum_{i\in I}trace_{F_{i}/F}(\tau_{i}(w_{i})w'_{i}x_{i}).$$
 On a les m\^emes propri\'et\'es que dans le cas $d$ pair. Pour un \'el\'ement semi-simple  fortement r\'eguli\`er $\tilde{x}\in \tilde{G}(F)$, son param\`etre $x_{D}$ n'est d\'etermin\'e que modulo le groupe des carr\'es $F^{\times,2}$.
 
 \bigskip
 
 {\bf Le cas du groupe unitaire.} Donnons-nous la collection d'objets suivante:
 
 $\bullet$ un ensemble fini $I$;
 
 $\bullet$ pour tout $i\in I$, une extension finie $F_{\pm i}$ de $F$; on pose $F_{i}=F_{\pm i}\otimes_{F}E$; c'est une $E$-alg\`ebre commutative; l'unique $F_{\pm i}$-automorphisme non trivial $\tau_{i}$ de $F_{i}$ est $id\otimes \tau_{E/F}$;
 
 $\bullet$ pour tout $i\in I$,  un \'el\'ement $c_{i}\in F_{i}^{\times}$ tel que $\tau_{i}(c_{i})=c_{i}$ et un \'el\'ement $x_{i}\in F_{i}^{\times}$ tel que $x_{i}\tau_{i}(x_{i})=1$.
 
 On pose $W=\sum_{i\in I}F_{i}$. C'est un espace vectoriel sur $E$ de dimension $\sum_{i\in I}[F_{i}:E]$. On suppose que cette dimension est \'egale \`a $d$. On munit $W$ de la forme hermitienne $q_{W}$ d\'efinie par
 $$(3) \qquad q_{W}(\sum_{i\in I}w_{i},\sum_{i\in I}w'_{i})=\sum_{i\in I}trace_{F_{i}/E}(\tau_{i}(w_{i})w'_{i}c_{i}).$$
 On suppose que $(W,q_{W})$ est $E$-isomorphe \`a $(V,q)$ et on fixe un tel isomorphisme. Modulo celui-ci, on d\'efinit un \'el\'ement $x\in G(F)$ par la formule (2). On a alors les m\^emes propri\'et\'es que dans le cas symplectique (quant \`a l'unicit\'e des param\`etres, c'est cette fois la classe d'isomorphie sur $E$ des couples $(F_{i},x_{i})$ qui est uniquement d\'etermin\'ee).
 
 \bigskip
 
 {\bf Le cas du changement de base du groupe unitaire.} Donnons-nous la collection d'objets suivante:
 
 $\bullet$ un ensemble fini $I$;
 
 $\bullet$ pour tout $i\in I$, une extension finie $F_{\pm i}$ de $F$; on pose $F_{i}=F_{\pm i}\otimes_{F}E$;
 
 $\bullet$ pour tout $i\in I$, un \'el\'ement $x_{i}\in F_{i}^{\times}$.
 
 On suppose $d=\sum_{i\in I}[F_{i}:E]$. On fixe un $E$-isomorphisme de $V$ sur $\sum_{i\in I}F_{i}$. Modulo cet isomorphisme, on d\'efinit $\tilde{x}\in \tilde{G}(F)$ par l'\'egalit\'e
 $$\tilde{x}(\sum_{i\in I}w_{i},\sum_{i\in I}w'_{i})=\sum_{i\in I}trace_{F_{i}/E}(\tau_{i}(w_{i})w'_{i}x_{i}).$$
 On a les m\^emes propri\'et\'es que dans le cas symplectique.  Pour une classe semi-simple fortement  r\'eguli\`ere $\tilde{x}\in \tilde{G}(F)$, ses param\`etres $I$ et $(F_{i})_{i\in I}$ sont essentiellement uniquement d\'etermin\'es (ce sont les classes de $E$-isomorphismes des $F_{i}$ qui comptent). Les $x_{i}$ sont d\'etermin\'es \`a multiplication pr\`es par le groupe $Norm_{F_{i}/F_{\pm i}}(F_{i}^{\times})$.

 \bigskip
 
 {\bf Remarque.} En admettant que l'intersection de l'ensemble des lecteurs de cet article et celui des lecteurs de [W] soit non vide, les \'el\'ements de cette intersection prendront garde au fait que les formules (1) et (3) ne co\"{\i}ncident pas avec celles de [W]. Dans cette r\'ef\'erence, on avait ajout\'e des coefficients $[F_{i}:F]^{-1}$ ou $[F_{i}:E]^{-1}$.
 
 \bigskip
 
 \subsection{Classes de conjugaison stable d'\'el\'ements semi-simples  suffisamment r\'eguliers}
 
 Deux \'el\'ements semi-simples suffisamment r\'eguliers sont stablement conjugu\'es si et seulement s'ils sont conjugu\'es par un \'el\'ement de $G(\bar{F})$. Le param\'etrage des classes de conjugaison stable de tels \'el\'ements se d\'eduit ais\'ement de celui du paragraphe pr\'ec\'edent. Dans les cas symplectique, sp\'ecial orthogonal ou unitaire, il suffit d'oublier les donn\'ees $(c_{i})_{i\in I}$ (dans le cas sp\'ecial orthogonal pair, on classifie ainsi des couples de classes de conjugaison stable). Dans les cas du groupe lin\'eaire tordu ou du changement de base du groupe unitaire, on remplace les classes $x_{i}Norm_{F_{i}/F_{\pm i}}(F_{i}^{\times})$ (ce sont ces classes qui param\'etraient les classes de conjugaison) par les classes $x_{i}F_{\pm i}^{\times}$. Dans le cas du groupe lin\'eaire tordu avec $d$ impair, on oublie de plus l'\'el\'ement $x_{D}$.
\bigskip

\subsection{Les formes quasi-d\'eploy\'ees}

Pour d\'efinir des facteurs de transfert, il importe de fixer un groupe quasi-d\'eploy\'e $\underline{G}$ et un torseur int\'erieur $\psi:G\to \underline{G}$. Que $\psi$ soit un torseur int\'erieur signifie que c'est un isomorphisme sur $\bar{F}$ et  qu'il existe une application $u:Gal(\bar{F}/F)\to \underline{G}(\bar{F})$ telle que $\sigma(\psi)\psi^{-1}(\underline{g})=u(\sigma)\underline{g}u(\sigma)^{-1}$ pour tout $\underline{g}\in \underline{G}(\bar{F})$. En g\'en\'eral,  l'application $u$, compos\'ee avec l'application naturelle de $\underline{G}(\bar{F})$ dans son groupe adjoint,  est un cocycle. Dans les cas qui nous int\'eressent, on peut effectuer les choix de sorte que $u$ elle-m\^eme soit un cocycle \`a valeurs dans $\underline{G}(\bar{F})$.  D'apr\`es une remarque de Kottwitz, cela permet de d\'efinir des facteurs de transfert pour le groupe $G$ de la m\^eme fa\c{c}on que si $G$ lui-m\^eme \'etait quasi-d\'eploy\'e. Pr\'ecis\'ement, il faut fixer le cocycle $u$ et le facteur de transfert d\'epend de $u$ et pas seulement de sa classe de cohomologie. Dans le cas symplectique, ou celui du groupe lin\'eaire tordu, ou celui du changement de base du groupe unitaire, $G$ est quasi-d\'eploy\'e, on choisit $\underline{G}=G$, $\psi$ est l'identit\'e et $u(\sigma)=1$ pour tout $\sigma\in Gal(\bar{F}/F)$. Pour unifier les notations, on pose dans ces cas $\underline{V}=V$ et, dans le cas symplectique, $\underline{q}=q$.

{\bf Le cas sp\'ecial orthogonal impair.} A \'equivalence pr\`es, il existe un unique espace $\underline{V}$  sur $F$, muni d'une forme quadratique non d\'eg\'en\'er\'ee $\underline{q}$, qui v\'erifie les conditions suivantes: la dimension de $\underline{V}$ sur $F$ est $d$; les discriminants de $q$ et $\underline{q}$ sont \'egaux dans $F^{\times}/F^{\times,2}$; le groupe sp\'ecial orthogonal $\underline{G}$ de $(\underline{V},\underline{q})$ est quasi-d\'eploy\'e. On introduit cet  espace quadratique. Les formes $q$ et $\underline{q}$ d\'efinissent par lin\'earit\'e des formes quadratiques sur $V\otimes_{F}\bar{F}$ et $\underline{V}\otimes_{F}\bar{F}$ (\`a valeurs dans $\bar{F}$). Ces deux formes sont isomorphes. On fixe un isomorphisme $\bar{F}$-lin\'eaire $\delta:V\otimes_{F}\bar{F}\to \underline{V}\otimes_{F}\bar{F}$ tel que
$\underline{q}(\delta(v),\delta(v'))=q(v,v')$ pour tous $v,v'\in V\otimes_{F}\bar{F}$. On d\'efinit un isomorphisme $\psi:G\to \underline{G}$ par $\psi(g)=\delta g\delta^{-1}$.  C'est un torseur int\'erieur. On  prend  pour fonction $u$  la fonction d\'efinie par $u(\sigma)=\sigma(\delta)\delta^{-1}$. Cette application $u$ est un cocycle \`a valeurs dans $\underline{G}(\bar{F})$.
\bigskip

{\bf Le cas sp\'ecial orthogonal pair.} Il est identique au pr\'ec\'edent, \`a ceci pr\`es que le couple $(\underline{V},\underline{q})$ n'est pas forc\'ement unique. On en choisit un et on construit un cocycle $u$ et un torseur int\'erieur $\psi:G\to \underline{G}$ comme ci-dessus. Remarquons que, m\^eme si $G$ est quasi-d\'eploy\'e, cette construction peut nous fournir un cocycle $u$ non trivial.

\bigskip

{\bf Le cas unitaire.} Il existe un  espace $\underline{V}$  sur $E$, muni d'une forme hermitienne non d\'eg\'en\'er\'ee $\underline{q}$, qui v\'erifie les conditions suivantes: la dimension de $\underline{V}$ sur $E$ est $d$;  le groupe  unitaire $\underline{G}$ de $(\underline{V},\underline{q})$ est quasi-d\'eploy\'e. Le couple $(\underline{V},\underline{q})$ est unique si $d$ est pair, il y en a deux si $d$ est impair. On choisit un tel  espace hermitien. Les formes $q$ et $\underline{q}$ d\'efinissent par lin\'earit\'e des formes $E\otimes_{F}\bar{F}$-hermitiennes sur $V\otimes_{F}\bar{F}$ et $\underline{V}\otimes_{F}\bar{F}$. Ces deux formes sont isomorphes. On fixe un isomorphisme $E\otimes_{F}\bar{F}$-lin\'eaire $\delta:V\otimes_{F}\bar{F}\to \underline{V}\otimes_{F}\bar{F}$ tel que
$\underline{q}(\delta(v),\delta(v'))=q(v,v')$ pour tous $v,v'\in V\otimes_{F}\bar{F}$. On d\'efinit un isomorphisme $\psi:G\to \underline{G}$ par $\psi(g)=\delta g\delta^{-1}$.  C'est un torseur int\'erieur. On  prend  pour fonction $u$  la fonction d\'efinie par $u(\sigma)=\sigma(\delta)\delta^{-1}$. Cette application $u$ est un cocycle \`a valeurs dans $\underline{G}(\bar{F})$.
\bigskip

\subsection{Epinglage}

Consid\'erons d'abord les cas symplectique, sp\'ecial orthogonal pair ou impair, ou unitaire. On fixe un sous-groupe de Borel $\underline{B}$ de $\underline{G}$ d\'efini sur $F$, un sous-tore maximal $\underline{T}$ de $\underline{B}$ et un \'epinglage invariant par $Gal(\bar{F}/F)$ relatif \`a la paire $(\underline{B},\underline{T})$. C'est-\`a-dire que, pour toute racine simple $\alpha$ de $\underline{T}$ dans l'alg\`ebre de Lie du radical unipotent de $\underline{B}$, on fixe un \'el\'ement non nul $E_{\alpha}$ du sous-espace radiciel correspondant, et on suppose que $E_{\sigma(\alpha)}=\sigma(E_{\alpha})$ pour toute racine simple $\alpha$ et tout $\sigma\in Gal(\bar{F}/F)$. Notons $N$ la somme des $E_{\alpha}$. C'est un \'el\'ement nilpotent r\'egulier de $\mathfrak{g}(F)$. Dans le cas symplectique ou sp\'ecial orthogonal impair, consid\'erons la forme bilin\'eaire $(v,v')\mapsto \underline{q}(v,N^{d-1}v')$ sur $\underline{V}$. Elle est quadratique de rang $1$. Elle est donc \'equivalente \`a la somme d'une forme nulle et d'une forme sur $F$ du type $(x,y)\mapsto \eta xy$, avec $\eta\in F^{\times}$. Cela d\'efinit un unique $\eta$ modulo le groupe des carr\'es $F^{\times,2}$. Dans le cas sp\'ecial orthogonal pair, on consid\`ere la forme  bilin\'eaire $(v,v')\mapsto \underline{q}(v,N^{d-2}v')$. La m\^eme construction s'applique. Dans le cas unitaire, on consid\`ere la forme sesquilin\'eaire  $(v,v')\mapsto \underline{q}(v,N^{d-1}v')$. Elle est hermitienne si $d$ est impair, antihermitienne si $d$ est pair. Elle est de rang $1$ et donc \'equivalente \`a la somme d'une forme nulle et d'une forme sur $E$ du type $(x,y)\mapsto \eta\tau_{E/F}(x)y$ avec $\eta\in E^{\times}$ (on a $\eta\in F^{\times}$ si $d$ est impair, $\tau_{E/F}(\eta)=-\eta$ si $d$ est pair). Cela d\'efinit un unique $\eta$ modulo le groupe des normes $Norm_{E/F}(E^{\times})$.

Consid\'erons maintenant les cas du groupe lin\'eaire tordu ou du changement de base du groupe unitaire. Modulo les choix d'une base de $V$ et d'un scalaire $\nu$, on a d\'efini en 1.2 un \'el\'ement $\tilde{\theta}$  de $\tilde{G}(F)$. Dans le cas du groupe lin\'eaire tordu, c'est une forme symplectique si $d$ est pair, quadratique si $d$ est impair. Dans le cas du changement de base du groupe unitaire, elle est proportionnelle \`a une forme hermitienne.  Notons $G_{\tilde{\theta}}$  la composante neutre de son groupe d'automorphismes. C'est un groupe  quasi-d\'eploy\'e. On fixe comme ci-dessus une "paire de Borel \'epingl\'ee" de ce groupe, invariante par $Gal(\bar{F}/F)$ et on en d\'eduit un \'el\'ement nilpotent $N\in \mathfrak{g}_{\tilde{\theta}}(F)$. On d\'efinit  $\eta$ comme pr\'ec\'edemment en consid\'erant la forme $(v,v')\mapsto  \tilde{\theta}(v,N^{d-1}v')$ sur $V$. Dans le cas du groupe lin\'eaire tordu, $\eta$ est un \'el\'ement de $F^{\times}$ bien d\'efini modulo $F^{\times,2}$. Dans le cas du changement de base du groupe unitaire, c'est un \'el\'ement de $E^{\times}$ bien d\'efini modulo $Norm_{E/F}(E^{\times})$ (ici, $\eta$ peut \^etre quelconque dans $E^{\times}$ puisqu'il d\'epend du scalaire $\nu$).

\bigskip

\subsection{$L$-groupes}

 Soit $N\geq1$ un entier. Le groupe $GL_{N}({\mathbb C})$ est le groupe des automorphismes de ${\mathbb C}^{N}$. On note $(\hat{e}_{k})_{k=1,...,N}$ la base standard de cet espace. On note $\hat{\theta}_{N}$ l'automorphisme de $GL_{N}({\mathbb C})$ d\'efini par la m\^eme formule qu'en 1.2, c'est-\`a-dire $\hat{\theta}_{N}(g)=J_{N}{^tg}^{-1}J_{N}^{-1}$. 
 On note $I^-_{N}$ l'\'el\'ement de $GL_{N}({\mathbb C})$ qui est antidiagonal et a pour coefficients antidiagonaux 
 $$(I^-_{N})_{k,N+1-k}=(-1)^{k}. $$
 Si $N$ est impair, on pose $I^+_{N}=I^-_{N}$. 
 Si $N$ est pair, on note $I^+_{N}$   l'\'el\'ement de $GL_{N}({\mathbb C})$ qui est antidiagonal et a pour coefficients antidiagonaux
  $$(I^+_{N})_{k,N+1-k}=\left\lbrace\begin{array}{cc}(-1)^{k},&\text{ si }k\leq N/2,\\ (-1)^{ k+1},&\text{ si }k>N/2.\\ \end{array}\right.$$
 On note $SO_{N}({\mathbb C})$ la composante neutre du sous-groupe des $g\in GL_{N}({\mathbb C})$ tels que $^tgI^+_{N}g=I^+_{N}$. Si $N$ est pair, on note $Sp_{N}({\mathbb C})$ le sous-groupe des  $g\in GL_{N}({\mathbb C})$ tels que $^tgI^-_{N}g=I^-_{N}$.

{\bf Cas symplectique.} Le $L$-groupe de $G$ est le produit direct $SO_{d+1}({\mathbb C})\times W_{F}$.

{\bf Cas sp\'ecial orthogonal impair.} Le $L$-groupe de $G$ est le produit direct $Sp_{d-1}({\mathbb C})\times W_{F}$.

{\bf Cas sp\'ecial orthogonal pair.}  Soit $\delta$ le discriminant de $q$. Si $\delta=1$, le $L$-groupe de $G$ est le produit direct $SO_{d}({\mathbb C})\times W_{F}$. Si $\delta\not=1$, on note $E=F(\sqrt{\delta})$. Le $L$-groupe de $G$ est le produit semi-direct $SO_{d}({\mathbb C})\rtimes W_{F}$. Un \'el\'ement de $W_{E}$ agit trivialement sur $SO_{d}({\mathbb C})$. Un \'el\'ement de $W_{F}\setminus W_{E}$ agit par conjugaison par  la matrice de permutation de ${\mathbb C}^{d}$ qui \'echange les vecteurs  $\hat{e}_{d/2}$ et $\hat{e}_{d/2+1}$ et fixe les autres vecteurs de base.

{\bf Cas du groupe lin\'eaire tordu.} Le $L$-groupe de $G$ est bien s\^ur $GL_{d}({\mathbb C})$. L'automorphisme "dual" de $\theta_{d}$ est $\hat{\theta}_{d}$.

{\bf Cas du groupe unitaire.} Le $L$-groupe de $G$ est le produit semi-direct $GL_{d}({\mathbb C})\rtimes W_{F}$. Un \'el\'ement de $W_{E}$ agit trivialement sur $GL_{d}({\mathbb C})$. Un \'el\'ement de $W_{F}\setminus W_{E}$ agit par l'automorphisme $\hat{\theta}_{d}$.

{\bf Cas du changement de base du groupe unitaire.} Le $L$-groupe de $G$ est le produit semi-direct $(GL_{d}({\mathbb C})\times GL_{d}({\mathbb C}))\rtimes W_{F}$. Un \'el\'ement de $W_{E}$ agit trivialement sur $GL_{d}({\mathbb C})\times GL_{d}({\mathbb C})$. Un \'el\'ement de $W_{F}\setminus W_{E}$ agit par $(g',g'')\mapsto (g'',g')$. L'automorphisme dual de $\theta_{d,E/F}$ est $(g',g'')\mapsto (\hat{\theta}_{d}(g''),\hat{\theta}_{d}(g'))$.

\bigskip

\subsection{Donn\'ees endoscopiques elliptiques}

On fixe une donn\'ee endoscopique elliptique pour $G$, ou pour $\tilde{G}$ dans le cas du groupe lin\'eaire tordu ou du changement de base du groupe unitaire. Une telle donn\'ee  est un triplet $(H,s,{^L{\xi}})$. Le premier terme $H$ est un groupe r\'eductif connexe d\'efini et quasi-d\'eploy\'e sur $F$, le deuxi\`eme est un \'el\'ement semi-simple de la composante complexe $\hat{G}$ du  $L$-groupe $^LG$ et le troisi\`eme est un plongement de $L$-groupes $^L\xi:{^LH}\to {^LG}$. Dans nos cas, $H$ est un produit $H^-\times H^+$, o\`u $H^-$ et $H^+$ sont chacun de l'un des cas symplectique, sp\'ecial orthogonal ou unitaire. On peut repr\'esenter $H^-$, resp. $H^+$, comme un sous- groupe du groupe d'automorphismes d'un couple $(V^-,q^-)$, resp. $(V^+,q^+)$ comme en 1.2. On note $d^-$, resp. $d^+$, la dimension de $V^-$, resp. $V^+$. L'\'el\'ement $s$ est toujours un \'el\'ement diagonal (dans la repr\'esentation du groupe $\hat{G}$ d\'ecrite dans le paragraphe pr\'ec\'edent), ou un produit $(s',s'')$ d'\'el\'ements diagonaux dans le cas du changement de base du groupe unitaire. Ces \'el\'ements diagonaux n'ont pour coefficients diagonaux  $s_{k}$ (ou $s'_{k}$ et $s''_{k}$) que des $\pm 1$. Pour d\'ecrire $^L\xi$, on repr\'esente les $L$-groupes de $H^-$, $H^+$ et $G$ comme en 1.7. En particulier, le groupe $\hat{H}^-$,  resp. $\hat{H}^+$, $\hat{G}$, est un sous-groupe du groupe des automorphismes d'un espace complexe dont on a fix\'e une base. On note $\hat{e}_{k}^-$, resp. $\hat{e}_{k}^+$, $\hat{e}_{k}$, les \'el\'ements de cette base. La restriction de $^L\xi$ \`a $W_{F}$ s'\'ecrit $w\mapsto( \rho(w),w)$, avec $\rho(w)\in \hat{G}$.

{\bf Cas symplectique.} Le groupe $H^-$ est sp\'ecial orthogonal pair.  On note $\delta^-$ le discriminant de $q^-$ et, si $\delta^-\not=1$, $E^-$ l'extension $F(\sqrt{\delta^-})$. La condition d'ellipticit\'e exclut le cas o\`u $d^-=2$ et $\delta^-=1$. Le groupe $H^+$ est symplectique. On a $d^-+d^+=d$. On a
$$s_{k}=\left\lbrace\begin{array}{cc}-1,&\text{ si }k=1,...,d^-/2\text { ou } k=d-d^-/2+1,...,d,\\ 1,&\text{ si }k=d^-/2+1,...,d-d^-/2.\\ \end{array}\right.$$
 La restriction de $^L\xi$ \`a $\hat{H}^-\times \hat{H}^+$ se d\'eduit de l'identification suivante des bases:
 $$\hat{e}_{k}^-\mapsto \left\lbrace\begin{array}{cc}\hat{e}_{k},&\text{ si }k=1,...,d^-/2,\\ \hat{e}_{k+d^+},&\text{ si }k=d^-/2+1,...,d^-;\\ \end{array}\right.$$
$$\hat{e}_{k}^+\mapsto \hat{e}_{k+d^-/2}.$$
Si $\delta^-=1$, on a $\rho(w)=1$ pour tout $w\in W_{F}$. Si $\delta^-\not=1$, on a $\rho(w)=1$ pour $w\in W_{E^-}$. Pour $w\in W_{F}\setminus W_{E^-}$,
$$\rho(w)\hat{e}_{k}=\left\lbrace\begin{array}{cc}\hat{e}_{k},&\text{ si }k=1,...,d^-/2-1\text{ ou }k=d-d^-/2+2,...,d,\\ \hat{e}_{d-d^-/2+1},&\text{ si }k=d^-/2,\\ -\hat{e}_{k},&\text{ si }k=d^-/2+1,...,d-d^-/2,\\ \hat{e}_{d^-/2},&\text{ si }k=d-d^-/2+1.\\ \end{array}\right.$$

{\bf Cas sp\'ecial orthogonal impair.} Les groupes $H^-$ et $H^+$ sont sp\'eciaux orthogonaux impairs. On a $d^-+d^+=d+1$. On a
$$s_{k}=\left\lbrace\begin{array}{cc}-1,&\text{ si }k=1,...,(d^--1)/2\text{ ou } k=d-(d^--1)/2,...,d-1,\\ 1,&\text{ si }k=(d^-+1)/2,...,d-(d^-+1)/2.\\ \end{array}\right.$$
 La restriction de $^L\xi$ \`a $\hat{H}^-\times \hat{H}^+$ se d\'eduit de l'identification suivante des bases:
 $$\hat{e}_{k}^-\mapsto \left\lbrace\begin{array}{cc}\hat{e}_{k},&\text{ si }k=1,...,(d^--1)/2,\\ \hat{e}_{k+d^+-1},&\text{ si }k=(d^-+1)/2,...,d^--1;\\ \end{array}\right.$$
$$\hat{e}_{k}^+\mapsto \hat{e}_{k+(d^--1)/2}.$$
On a $\rho(w)=1$ pour tout $w\in W_{F}$.

{\bf Cas sp\'ecial orthogonal pair.} Les groupes $H^-$ et $H^+$ sont sp\'eciaux orthogonaux pairs. On a $d^-+d^+=d$. On note $\delta$, $\delta^-$ et $\delta^+$ les discriminants de $q$, $q^-$ et $q^+$ et $E$, $E^-$ et $E^+$ les extensions associ\'ees, quand elles existent.  La condition d'ellipticit\'e exclut les cas o\`u l'un au moins des couples$(d^-,\delta^-)$ ou $(d^+,\delta^+)$ est \'egal \`a $(2,1)$. On a $\delta^-\delta^+=\delta$. On a
$$s_{k}=\left\lbrace\begin{array}{cc}-1,&\text{ si }k=1,...,d^-/2\text{ ou } k=d-d^-/2+1,...,d,\\ 1,&\text{ si }k=d^-/2+1,...,d-d^-/2.\\ \end{array}\right.$$
 La restriction de $^L\xi$ \`a $\hat{H}^-\times \hat{H}^+$ se d\'eduit de l'identification suivante des bases:
 $$\hat{e}_{k}^-\mapsto \left\lbrace\begin{array}{cc}\hat{e}_{k},&\text{ si }k=1,...,d^-/2,\\ \hat{e}_{k+d^+},&\text{ si }k=d^-/2+1,...,d^-;\\ \end{array}\right.$$
$$\hat{e}_{k}^+\mapsto \hat{e}_{k+d^-/2}.$$
Si $d^+\not=0$, notons $S\in GL_{d}({\mathbb C})$ la matrice de permutation qui \'echange les vecteurs $\hat{e}_{d^-/2}$ et $\hat{e}_{d-d^-/2+1}$, ainsi que les vecteurs $\hat{e}_{d/2}$ et $\hat{e}_{d/2+1}$, et qui fixe les autres vecteurs de base. Si $d^+=0$, $S=1$. Si $\delta^-=1$, $\rho(w)=1$ pour tout $w\in W_{F}$. Si $\delta^-\not=1$, $\rho(w)=1$ pour $w\in W_{E^-}$ et $\rho(w)=S$ pour $w\in W_{F}\setminus W_{E^-}$.

{\bf Cas du groupe lin\'eaire tordu avec $d$ pair.} Le groupe $H^-$ est sp\'ecial orthogonal pair. On note
 $\delta^-$ le discriminant de $q^-$ et, si $\delta^-\not=1$, $E^-$ l'extension $F(\sqrt{\delta^-})$. La condition d'ellipticit\'e exclut le cas o\`u $d^-=2$ et $\delta^-=1$. Le groupe $H^+$ est sp\'ecial orthogonal impair. On a $d^-+d^+=d+1$. On a
 $$s_{k}=\left\lbrace\begin{array}{cc}-1,&\text{ si }k=1,...,d^-/2 ,\\ 1,&\text{ si }k=d^-/2+1,...,d.\\ \end{array}\right.$$
  La restriction de $^L\xi$ \`a $\hat{H}^-\times \hat{H}^+$ se d\'eduit de l'identification suivante des bases:
 $$\hat{e}_{k}^-\mapsto \left\lbrace\begin{array}{cc}\hat{e}_{k},&\text{ si }k=1,...,d^-/2,\\ \hat{e}_{k+d^+-1},&\text{ si }k=d^-/2+1,...,d^-;\\ \end{array}\right.$$
$$\hat{e}_{k}^+\mapsto \hat{e}_{k+d^-/2}.$$
Si $\delta^-=1$, on a $\rho(w)=1$ pour tout $w\in W_{F}$. Si $\delta^-\not=1$, on a $\rho(w)=1$ pour tout $w\in W_{E^-}$. Pour $w\in W_{F}\setminus W_{E^-}$, $\rho(w)$ est la matrice de permutation qui \'echange les vecteurs $\hat{e}_{d^-/2}$ et $\hat{e}_{d-d^-/2+1}$ et qui fixe tout autre vecteur de base.

{\bf Cas du groupe lin\'eaire tordu avec $d$ impair.} Le groupe $H^-$ est sp\'ecial orthogonal impair. Le groupe $H^+$ est symplectique. On a $d^-+d^+=d$.  On a
 $$s_{k}=\left\lbrace\begin{array}{cc}-1,&\text{ si }k=1,...,(d^--1)/2 ,\\ 1,&\text{ si }k=(d^-+1)/2,...,d.\\ \end{array}\right.$$
  La restriction de $^L\xi$ \`a $\hat{H}^-\times \hat{H}^+$ se d\'eduit de l'identification suivante des bases:
 $$\hat{e}_{k}^-\mapsto \left\lbrace\begin{array}{cc}\hat{e}_{k},&\text{ si }k=1,...,(d^--1)/2,\\ \hat{e}_{k+d^++1},&\text{ si }k=(d^-+1)/2,...,d^--1;\\ \end{array}\right.$$
$$\hat{e}_{k}^+\mapsto \hat{e}_{k+(d^--1)/2}.$$
On choisit un caract\`ere $\chi$ de $F^{\times}$, d'ordre au plus $2$, que l'on identifie \`a un caract\`ere de $W_{F}$. Pour tout $w\in W_{F}$, $\rho(w)$ est la matrice diagonale de coefficients
$$\rho(w)_{k}=\left\lbrace\begin{array}{cc}1,&\text{ si }k=1,...,(d^--1)/2\text{ ou }k=d-(d^--1)/2+1,...,d,\\ \chi(w),&\text{ si }k=(d^-+1)/2,...,d-(d^--1)/2.\\ \end{array}\right.$$

{\bf Cas  unitaire.} Les groupes $H^-$ et $H^+$ sont unitaires, relatifs \`a la m\^eme extension $E$ que $G$. On a
$$s_{k}=\left\lbrace\begin{array}{cc}-1,&\text{ si }k=1,...,d^-\\ 1,&\text{ si }k=d^-+1,...,d.\\ \end{array}\right.$$
La restriction de $^L\xi$ \`a $\hat{H}^-\times \hat{H}^+$ est
$$(g^-,g^+)\mapsto \left(\begin{array}{cc}g^-&0\\ 0&g^+\\ \end{array}\right).$$
On choisit deux caract\`eres $\mu^-$ et $\mu^+$ de $E^{\times}$. On suppose que la restriction de $\mu^-$, resp. $\mu^+$, \`a $F^{\times}$ co\"{\i}ncide avec $sgn_{E/F}^{d^+}$, resp. $sgn_{E/F}^{d^-}$.  On identifie $\mu^-$ et $\mu^+$ \`a des caract\`eres de $W_{E}$. On fixe de plus deux nombres complexes non nuls $z^-$ et $z^+$. Pour $w\in W_{E}$, $\rho(w)=\boldsymbol{\mu}(w)$, o\`u $\boldsymbol{\mu}(w)$ est la matrice diagonale de coefficients diagonaux
$$ \boldsymbol{\mu}(w)_{k}=\left\lbrace\begin{array}{cc}\mu^-(w),&\text{ si } k=1,...,d^-,\\ \mu^+(w),&\text{ si }k=d^-+1,...,d.\\ \end{array}\right.$$
Pour $w\in W_{F}\setminus W_{E}$, on a $\rho(w)=Z$, o\`u $Z\in GL_{d}({\mathbb C})$ est d\'efinie par
$$Z\hat{e}_{k}=\left\lbrace\begin{array}{cc}z^+\hat{e}_{k+d^-},&\text{ si }k=1,...,d^+,\\ z^-\hat{e}_{k-d^+},&\text{ si }k=d^++1,...,d.\\ \end{array}\right.$$

{\bf Cas du changement de base du groupe unitaire.} Les groupes $H^-$ et $H^+$ sont unitaires, relatifs \`a la m\^eme extension $E$ que $G$. On a $s=(s',s'')$, avec $s''=1$ et
$$s'_{k}=\left\lbrace\begin{array}{cc}-1,&\text{ si }k=1,...,d^-\\ 1,&\text{ si }k=d^-+1,...,d.\\ \end{array}\right.$$ 
La restriction de $^L\xi$ \`a $\hat{H}^-\times \hat{H}^+$ est $(g^-,g^+)\mapsto (g',g'')$, o\`u
$$g'= \left(\begin{array}{cc}g^-&0\\ 0&g^+\\ \end{array}\right), \,\, g''=\hat{\theta}_{d}(g').$$
On choisit  deux caract\`eres $\mu^-$ et $\mu^+$ de $E^{\times}$. On suppose que la restriction de $\mu^-$, resp. $\mu^+$, \`a $F^{\times}$ co\"{\i}ncide avec $sgn_{E/F}^{d^++1}$, resp. $sgn_{E/F}^{d^-}$.  On identifie $\mu^-$ et $\mu^+$ \`a des caract\`eres de $W_{E}$. On fixe de plus deux nombres complexes non nuls $z^-$ et $z^+$. On d\'efinit $Z$ et $\boldsymbol{\mu}(w)$ comme dans le cas unitaire. Pour $w\in W_{E}$, $\rho(w)=(\boldsymbol{\mu}(w),\hat{\theta}_{d}(\boldsymbol{\mu}(w)))$. Pour $w\in W_{F}\setminus W_{E}$, $\rho(w)=(Z, \hat{\theta}_{d}(Z)s')$.

\bigskip

On a effectu\'e divers choix qui sont n\'ecessaires pour d\'efinir les facteurs de transfert (certains sont toutefois inessentiels, par exemple ceux des \'el\'ements $z^+$ et $z^-$ dans les cas unitaire ou du changement de base du groupe unitaire). Toutefois, la plupart des choix n'affectent pas la classe d'\'equivalence de la donn\'ee $(H,s,{^L{\xi}})$. Indiquons quels sont les objets qui d\'eterminent cette classe.

{\bf Cas symplectique.} Le couple $(d^-,d^+)$ et le discriminant $\delta^-$.

{\bf Cas sp\'ecial orthogonal impair.} Le couple $(d^-,d^+)$ \`a l'ordre pr\`es, c'est-\`a-dire que $(d^+,d^-)$ est \'equivalent \`a $(d^-,d^+)$.

{\bf Cas sp\'ecial orthogonal pair.} Le quadruplet $(d^-,\delta^-,d^+,\delta^+)$, \`a la permutation suivante  pr\`es: $(d^+,\delta^+,d^-,\delta^-)$ est \'equivalent \`a $(d^-,\delta^-,d^+,\delta^+)$.

{\bf Cas du groupe lin\'eaire tordu, avec $d$ pair.} Le triplet $(d^-,\delta^-,d^+)$.

{\bf Cas du groupe lin\'eaire tordu, avec $d$ impair.} Le triplet $(d^-,d^+,\chi)$.

{\bf Cas du groupe unitaire.} Le couple $(d^-,d^+)$ \`a l'ordre pr\`es.

{\bf Cas du changement de base du groupe unitaire.} Le couple $(d^-,d^+)$. 

\bigskip

\subsection{Correspondance de classes de conjugaison stable}

On fixe un \'el\'ement $x\in G(F)$, ou $\tilde{x}\in \tilde{G}(F)$ dans les cas du groupe lin\'eaire tordu ou du changement de base du groupe unitaire, et un \'el\'ement $y=(y^-,y^+)\in H(F)$. On les suppose semi-simples et suffisamment r\'eguliers.  Introduisons des donn\'ees qui param\`etrent les classes de conjugaison stable de $y^+$ et $y^-$, cf. 1.4. On les note   $(I^-, (F_{\pm i})_{i\in I^-},(F_{i})_{i\in I^-}, (y_{i})_{i\in I^-})$ et $(I^+,(F_{\pm i})_{i\in I^+},(F_{i})_{i\in I^+},(y_{i})_{i\in I^+})$, en supposant que $I^-$ et $I^+$ sont des ensembles disjoints.  Supposons que les classes de conjugaison stable de $x$, resp. $\tilde{x}$, et de $y$ se correspondent. Alors on peut param\'etrer la classe de conjugaison stable de $x$, resp. $\tilde{x}$, par les donn\'ees suivantes:

$\bullet$ l'ensemble $I=I^-\cup I^+$;

$\bullet$ pour $i\in I$, le corps $F_{\pm i}$ et la $F_{\pm i}$-alg\`ebre $F_{i}$ qui figurent dans les donn\'ees associ\'ees \`a $y^-$ ou $y^+$;

$\bullet$ pour $i\in I$, un \'el\'ement $x_{i}\in F_{i}^{\times}$ qui v\'erifie les conditions suivantes:

- dans les cas symplectique, sp\'ecial orthogonal ou unitaire, $x_{i}=y_{i}$;

- dans les cas du groupe lin\'eaire tordu ou du changement de base du groupe unitaire, $x_{i}\tau_{i}(x_{i})^{-1}=(-1)^{d+1}y_{i}\nu/\tau_{i}(\nu)$ (rappelons que $\nu$ figure dans la d\'efinition de l'\'el\'ement $\tilde{\theta}$) .

\bigskip

Hormis le cas sp\'ecial orthogonal pair, la r\'eciproque est vraie: si la classe de conjugaison stable de $x$, resp. $\tilde{x}$, est param\'etr\'ee par les donn\'ees ci-dessus, les classes de conjugaison stable de $x$ et de $y$ se correspondent. Dans le cas sp\'ecial orthogonal pair, on  introduit un \'el\'ement $x'$ qui est conjugu\'e \`a $x$ par un \'el\'ement du groupe orthogonal de d\'eterminant $-1$. On se rappelle que $x'$ n'est pas stablement conjugu\'e \`a $x$, mais les classes de conjugaison stable de $x$ et $x'$ sont param\'etr\'ees pas les m\^emes donn\'ees. Alors, soit les classes de conjugaison stable de $x$ et $y$ se correspondent, soit ce sont celles de $x'$ et $y$, ces deux cas \'etant exclusifs l'un de l'autre.

\bigskip

\subsection{Le facteur de transfert}

On a fix\'e un torseur int\'erieur $\psi:G\mapsto \underline{G}$ et un cocycle $u$ \`a valeurs dans $G(\bar{F})$. On a fix\'e un \'epinglage de $\underline{G}$, ou d'un certain sous-groupe dans les cas du groupe lin\'eaire tordu ou du changement de base du groupe unitaire. On a fix\'e une donn\'ee endoscopique $(H,s,{^L\xi})$ de $G$ ou $\tilde{G}$. C'est tout ce qu'il nous faut pour d\'efinir un facteur de transfert $\Delta_{H,G}$, ou $\Delta_{H,\tilde{G}}$ ([LS], [KS]). Nous supprimons de ce facteur les termes $\Delta_{IV}$. On  consid\`ere le facteur de transfert comme une fonction d\'efinie sur les couples $(y,x) \in H(F)\times G(F)$, resp. $(y,\tilde{x})\in H(F)\times \tilde{G}(F)$, form\'es d'\'el\'ements semi-simples suffisamment r\'eguliers et dont les classes de conjugaison stable se correspondent. 

Consid\'erons un tel couple. On \'ecrit $y=(y^-,y^+)$ et on param\`etre les classes de conjugaison stable de $y^-$ et $y^+$ comme dans le paragraphe pr\'ec\'edent. Comme en 1.3, la classe de conjugaison de $x$ est param\'etr\'ee par un ensemble $I$ et des familles $(F_{\pm i})_{i\in I}$, $(F_{i})_{i\in I}$, $(x_{i})_{i\in I}$ et, dans certains cas, une famille $(c_{i})_{i\in I}$ ou un \'el\'ement $x_{D}$. D'apr\`es le paragraphe pr\'ec\'edent, on peut supposer que $I=I^-\cup I^+$ et que les trois premi\`eres familles v\'erifient les conditions d\'ecrites dans ce paragraphe. Consid\'erons d'abord les cas symplectique ou sp\'ecial orthogonal ou du groupe lin\'eaire tordu. Pour tout $i\in I$, notons $\Phi_{i}$ l'ensemble des homomorphismes de $F$-alg\`ebres de $F_{i}$ dans $\bar{F}$. On d\'efinit le polyn\^ome 
$$P_{I}(T)=\prod_{i\in I}\prod_{\phi\in \Phi_{i}}(T-\phi(y_{i})).$$
On d\'efinit de fa\c{c}on similaire les polyn\^omes $P_{I^+}(T)$ et $P_{I^-}(T)$. On note $P_{I}'(T)$ le polyn\^ome d\'eriv\'e de $P_{I}(T)$. Consid\'erons maintenant les cas unitaire ou du changement de base du groupe unitaire. Pour tout $i\in I$, on note $\Phi_{i,E}$ l'ensemble des homomorphismes de $E$-alg\`ebres de $F_{i}$ dans $\bar{F}$. On d\'efinit le polyn\^ome $P_{I,E}(T)$ en rempla\c{c}ant $\Phi_{i}$ par $\Phi_{i,E}$ dans la formule ci-dessus.

Notons $I^*$ le sous-ensemble des $i\in I$ tels que $F_{i}$ soit un corps. D\'efinissons $I^{-*}$ et $I^{+*}$ de fa\c{c}on similaire. Soit $i\in I^{-*}$. On d\'efinit $C_{i}\in F_{i}^{\times}$ par les formules suivantes.

{\bf Cas symplectique.} $C_{i}=-\eta c_{i}P'_{I}(y_{i})P_{I}(-1)y_{i}^{1-d/2}$.

{\bf Cas sp\'ecial orthogonal impair .} 
$$C_{i}=-2\eta c_{i}P'_{I}(y_{i})P_{I}(-1)y_{i}^{(3-d)/2}(1+y_{i})(y_{i}-1)^{-1}.$$

{\bf Cas sp\'ecial orthogonal pair.}
$$C_{i}=2\eta c_{i}P'_{I}(y_{i})P_{I}(-1)y_{i}^{1-d/2}(1+y_{i})(y_{i}-1)^{-1}.$$

{\bf Cas du groupe lin\'eaire tordu, avec $d$ pair.}
$$C_{i}=\eta x_{i}^{-1}P'_{I}(y_{i})P_{I}(-1)y_{i}^{1-d/2}(1+y_{i}).$$

{\bf Cas du groupe lin\'eaire tordu, avec $d$ impair.} 
$$C_{i}=x_{D}x_{i}^{-1}P'_{I}(y_{i})P_{I}(1)y_{i}^{(3-d)/2}(y_{i}-1).$$

{\bf Cas du groupe unitaire, avec $d$ pair.} $C_{i}=-\eta c_{i}P'_{I,E}(y_{i})P_{I,E}(-1)^{-1}y_{i}^{1-d/2}$.

{\bf Cas du groupe unitaire, avec $d$ impair.}
$$C_{i}=-\eta c_{i}P'_{I,E}(y_{i})P_{I,E}(-1)^{-1}y_{i}^{(1-d)/2}(1+y_{i}).$$

{\bf Cas du changement de base du groupe unitaire, avec $d$ pair.}
$$C_{i}=-\eta x_{i}^{-1}P'_{I,E}(y_{i})P_{I,E}(-1)^{-1}y_{i}^{1-d/2}(1+y_{i}).$$

{\bf Cas du changement de base du groupe unitaire, avec $d$ impair.}
$$C_{i}=-\eta x_{i}^{-1}P'_{I,E}(y_{i})P_{I,E}(-1)^{-1}y_{i}^{(3-d)/2}.$$

Dans tous les cas, on v\'erifie que $C_{i}$ appartient \`a $F_{\pm i}^{\times}$.

\ass{Proposition}{On a les \'egalit\'es:

- dans les cas symplectique ou sp\'ecial orthogonal, resp. dans le cas du groupe lin\'eaire tordu, avec $d$ pair,
$$\left.\begin{array}{cc}&\Delta_{H,G}(y,x)\\ \text{ resp. }& \Delta_{H,\tilde{G}}(y,\tilde{x})\\ \end{array} \right\rbrace=\prod_{i\in I^{-*}}sgn_{F_{i}/F_{\pm i}}(C_{i});$$

- dans le cas du groupe lin\'eaire tordu, avec $d$ impair,
$$ \Delta_{H,\tilde{G}}(y,\tilde{x})=\chi(\eta x_{D}P_{I}(1)P_{I^-}(-1))\prod_{i\in I^{-*}}sgn_{F_{i}/F_{\pm i}}(C_{i});$$

- dans le cas unitaire, resp. dans le cas du changement de base du groupe unitaire,
$$\left.\begin{array}{cc}&\Delta_{H,G}(y,x)\\ \text{ resp. }& \Delta_{H,\tilde{G}}(y,\tilde{x})\\ \end{array} \right\rbrace=\mu^-(P_{I^-}(0)P_{I^-}(-1)^{-1})\mu^+(P_{I^+}(0)P_{I^+}(-1)^{-1})\prod_{i\in I^{-*}}sgn_{F_{i}/F_{\pm i}}(C_{i}).$$}

\bigskip

\subsection{Compl\'ements}

Dans les cas sp\'ecial orthogonal ou unitaire, on a remarqu\'e que la classe d'\'equivalence des donn\'ees $(H,s,{^L\xi})$ etait invariante par la permutation de $H^-$ et $H^+$. Cela signifie que, dans la proposition pr\'ec\'edente, on peut remplacer l'ensemble d'indices $I^{-*}$ par $I^{+*}$. On obtient encore un facteur de transfert. Si le cocycle $u$ est trivial, les deux facteurs ainsi obtenus sont \'egaux. Dans le cas o\`u $F$ est non archim\'edien,  il n'y a que deux classes de cocycles possibles. Alors, si $u$ est non trivial, les deux facteurs de transfert sont oppos\'es l'un de l'autre. 

Consid\'erons le cas du groupe lin\'eaire tordu, avec $d$ pair et $F$ non archim\'edien. Supposons $d^-=d$ et $d^+=1$. On a $H^+=\{1\}$. Le groupe $H^-$ est quasi-d\'eploy\'e. On en choisit une paire de Borel \'epingl\'ee invariante par $Gal(\bar{F}/F)$ et on en d\'eduit un invariant $\eta^-$ comme en 1.6. Supposons $\eta=-\eta^-$.  Consid\'erons un couple $(y,\tilde{x})$ comme dans le paragraphe pr\'ec\'edent. D\'efinissons une forme quadratique $q_{\tilde{x}}$ sur $V$ par 
$$q_{\tilde{x}}(v,v')=\tilde{x}(v,v')+\tilde{x}(v',v).$$
Parce que la classe de conjugaison stable de $\tilde{x}$ correspond \`a une classe de conjugaison stable de $H(F)$, la forme $q_{\tilde{x}}$ a m\^eme discriminant que $q^-$. On a les \'egalit\'es
$$\Delta_{H,\tilde{G}}(y,\tilde{x})=\left\lbrace\begin{array}{cc}1,&\text{ si } (V,q_{\tilde{x}})\text{ et }(V^-,q^-)\text{ sont isomorphes ,}\\ -1,&\text{ sinon. }\\ \end{array}\right.$$

\bigskip

\section{Preuve dans le cas du groupe lin\'eaire tordu, avec $d$ impair}

\bigskip

\subsection{Changement de $\tilde{x}$ dans sa classe de conjugaison}

On consid\`ere les donn\'ees de 1.10 dans le cas du groupe lin\'eaire tordu, avec $d$ impair.  L'\'el\'ement $\tilde{\theta}$ est une forme quadratique sur $V$.  On fixe un \'el\'ement $c_{D}\in F^{\times}$ et, pour tout $i\in I$, on fixe $c_{i}\in F_{\pm i}^{\times}$. On pose $D=F$. On suppose qu'il existe un isomorphisme de $V$ sur    $D\oplus(\oplus_{i\in I}F_{i})$, de sorte que, modulo cet isomorphisme, $\tilde{\theta}$ s'\'ecrive 
$$\tilde{\theta}(v_{D}+\sum_{i\in I}v_{i},v'_{D}+\sum_{i\in I}v'_{i})=c_{D}v_{D}v'_{D}+\sum_{i\in I}trace_{F_{i}/F}(\tau_{i}(v_{i})v'_{i}c_{i}).$$
De tels choix sont possibles. On fixe un tel isomorphisme. On note $G_{\tilde{\theta}}$ le groupe sp\'ecial orthogonal de $\tilde{\theta}$. Notons $T^{\flat}$ le sous-tore de $G_{\tilde{\theta}}$ form\'e des \'el\'ements qui respectent chaque sous-espace $F_{i}$ de $V$. On a $T^{\flat}(F)=\prod_{i\in I}F_{i}^1$, o\`u $F_{i}^1=\{\lambda\in F_{i}^{\times}; \lambda\tau_{i}(\lambda)=1\}$. Notons $T^{\diamond}$ le commutant de $T^{\flat}$ dans $G$. C'est aussi le sous-groupe des \'el\'ements de $G$ qui respectent chaque sous-espace $F_{i}$ ainsi que la droite $D$. On a $T^{\diamond}(F)=F^{\times}\times \prod_{i\in I}F_{i}^{\times}$.   Posons $t^{\diamond}_{x,D}=c_{D}x_{D}^{-1}$ et, pour $i\in I$, $t^{\diamond}_{x,i}=c_{i}\tau_{i}(x_{i})^{-1}$. D\'efinissons l'\'el\'ement $t^{\diamond}_{\tilde{x}}=(t^{\diamond}_{\tilde{x},D},(t^{\diamond}_{\tilde{x},i})_{i\in I}))$ de $T(F)$ et consid\'erons l'\'el\'ement $t^{\diamond}_{\tilde{x}}\tilde{\theta}$. Par d\'efinition de ce dernier produit, on a 
$(t^{\diamond}_{\tilde{x}}\tilde{\theta})(v,v')=\tilde{\theta}((t^{\diamond}_{\tilde{x}})^{-1}v,v')$ pour tout $v,v'\in V$. Plus explicitement
$$(t^{\diamond}_{\tilde{x}}\tilde{\theta})(v_{D}+\sum_{i\in I}v_{i},v'_{D}+\sum_{i\in I}v'_{i})=c_{D}(t^{\diamond}_{\tilde{x},D})^{-1}v_{D}v'_{D}+\sum_{i\in I}trace_{F_{i}/F}(\tau_{i}(v_{i})v'_{i}c_{i}\tau_{i}(t^{\diamond}_{\tilde{x},i})^{-1})$$
$$=x_{D}v_{D}v'_{D}+\sum_{i\in I}trace_{F_{i}/F}(\tau_{i}(v_{i})v'_{i} x_{i}).$$
On voit ainsi que $t^{\diamond}_{\tilde{x}}\tilde{\theta}$ est conjugu\'e \`a $\tilde{x}$.  Cela nous permet de supposer d\'esormais $\tilde{x}=t^{\diamond}_{\tilde{x}}\tilde{\theta}$.

\bigskip

\subsection{Syst\`emes de racines}

On pose simplement $\theta=\theta_{d}$, cf. 1.2. Modulo le choix d'une base de $V$, on identifie $G$ \`a $GL_{d}$ et $\tilde{G}$ \`a $GL_{d}\theta$ comme on l'a expliqu\'e dans ce paragraphe. Notons $T$ le sous-tore diagonal de $G$ et, pour $t\in T$, notons $(t_{k})_{k=1,...,d}$ ses coefficients diagonaux.    
 Notons $\Sigma$ l'ensemble des racines de $T$ dans $G$. Il s'identifie \`a l'ensemble des couples $(k,l)\in \{1,...,d\}^2$ tels que $k\not=l$: pour un tel couple $(k,l)$, la racine correspondante $\alpha_{k,l}$ 
est donn\'ee par $\alpha_{k,l}(t)=t_{k}t_{l}^{-1}$. L'automorphisme $\theta$ agit sur $\Sigma$. On a $\theta(\alpha_{k,l})=\alpha_{d+1-l,d+1,k}$. Pour $\alpha\in \Sigma$, on pose $N\alpha=\alpha+\theta(\alpha)$ si $\alpha\not=\theta(\alpha)$, $N\alpha=\alpha$ si $\alpha=\theta(\alpha)$. Kottwitz et Shelstad distinguent trois types de racines. On traduit explicitement leur description ainsi: une racine $\alpha_{k,l}$ est de type $R_{1}$ si $k\not=d+1-l$ et $k$ et $l$ sont tous deux diff\'erents de $(d+1)/2$, de type $R_{2}$ si $k$ ou $l$ est \'egal \`a $(d+1)/2$, de type $R_{3}$ si $k=d+1-l$.

Notons $T_{\theta}$ la composante neutre du groupe des points fixes de $\theta$ dans $T$. Un \'el\'ement $t\in T$ appartient \`a $T_{\theta}$ si et seulement si $t_{(d+1)/2}=1$ et $t_{d+1-k}=t_{k}^{-1}$ pour $k\not=(d+1)/2$. Le groupe $T_{\theta}$ est un sous-tore maximal de $G_{\tilde{\theta}}$. Notons $\Sigma_{\theta}$ l'ensemble des racines de $T_{\theta}$ dans $G_{\tilde{\theta}}$. Pour $\alpha\in \Sigma$, notons $\alpha_{res}$ la restriction de $\alpha$ \`a $T_{\theta}$.  Si $\alpha$ est de type $R_{1}$ ou $R_{2}$, $\alpha_{res}$ appartient \`a $\Sigma_{\theta}$. Si $\alpha$ est de type $R_{3}$, $\alpha_{res}/2$ appartient \`a $\Sigma_{\theta}$ (en adoptant une notation additive usuelle).

  Posons $\hat{G}=GL_{d}({\mathbb C})$, $\hat{\theta}=\hat{\theta}_{d}$, notons $\hat{T}$ le sous-tore diagonal de $\hat{G}$ et $ \hat{G}_{\hat{\theta}}$, resp. $\hat{T}_{\hat{\theta}}$, la composante neutre du sous-groupe des points fixes de $\hat{\theta}$ dans $\hat{G}$, resp. $\hat{T}$. On note $\check{\Sigma}$ l'ensemble des racines de $\hat{T}$ dans $\hat{G}$, qui s'identifie au m\^eme ensemble de couples que pr\'ec\'edemment. On note $ \check{\alpha}_{k,l}$ la racine associ\'ee \`a $(k,l)$. On note $\check{\Sigma}_{\hat{\theta}}$ l'ensemble des racines de $\hat{T}_{\hat{\theta}}$ dans $\hat{G}_{\hat{\theta}}$. Pour $\check{\alpha}\in \check{\Sigma}$, on note $\check{\alpha}_{res}$ sa restriction \`a $\hat{T}_{\hat{\theta}}$. Cette op\'eration de restriction v\'erifie des propri\'et\'es similaires \`a celles d\'ecrites ci-dessus. On a plong\'e $\hat{H}$ dans $\hat{G}$. Par ce plongement, $\hat{T}_{\hat{\theta}}$ s'identifie \`a un sous-tore maximal de $\hat{H}$. On note $\check{\Sigma}_{H}$ l'ensemble des racines de $\hat{T}_{\hat{\theta}}$ dans $\hat{H}$. Cet ensemble est r\'eunion de deux sous-ensembles \'evidents $\check{\Sigma}_{H^-}$ et $\check{\Sigma}_{H^+}$.
  
  L'ensemble $\check{\Sigma}$ s'identifie naturellement \`a l'ensemble de coracines associ\'e au syst\`eme de racines $\Sigma$. La bijection est \'evidemment $\alpha_{k,l}\mapsto \check{\alpha}_{k,l}$.  
    
  {\bf Attention:} l'ensemble $\check{\Sigma}_{\hat{\theta}}$ ne s'identifie pas \`a l'ensemble de coracines associ\'e \`a $\Sigma_{\theta}$.  Ces ensembles de racines sont tous deux de type $B_{d}$. Il y a toutefois une bijection naturelle entre $\Sigma_{\theta}$ et $\check{\Sigma}_{\hat{\theta}}$: pour $\beta\in \Sigma_{\theta}$, on choisit $\alpha\in \Sigma$ tel que $\beta=\alpha_{res}$; on associe \`a $\beta$ l'\'el\'ement $(\check{\alpha})_{res}$ de $\check{\Sigma}_{\hat{\theta}}$. Cela ne d\'epend pas du choix de $\alpha$.
  
  \bigskip
  
  \subsection{Description galoisienne du tore}
  
  Pour une extension finie $F'$ de $F$, on note  $\Phi_{F'}$ l'ensemble des homomorphismes de $F$-alg\`ebres de $F'$ dans $\bar{F}$. On a dit que l'on consid\'erait $F'$ comme un sous-corps de $\bar{F}$. Il revient au m\^eme de dire que l'on fixe un \'el\'ement privil\'egi\'e $\phi_{F'}\in \Phi_{F'}$, \`a savoir l'identit\'e.
    
      Pour $i\in I^*$, on pose $\Phi_{i}=\Phi_{F_{i}}$ et $\phi_{i}=\phi_{F_{i}}$. Pour $i\in I\setminus I^*$,  on note $\Phi_{i}$ l'ensemble des homomorphismes de $F$-alg\`ebres de $F_{i}$ dans $F$. Il y a deux homomorphismes de $F_{\pm i}$-alg\`ebres de $F_{i}$ dans $F_{\pm i}$, notons-les $\psi_{i}^1$ et $\psi_{i}^2$. Alors $\Phi_{i}=\{\phi\circ \psi_{i}^1; \phi\in\Phi_{F_{\pm i}}\}\sqcup\{\phi\circ \psi_{i}^2; \phi\in\Phi_{F_{\pm i}}\}$. On note simplement $\Phi_{i}^1$ et $\Phi_{i}^2$ les deux termes de cette d\'ecomposition et on pose $\phi_{i}^1=\phi_{F_{\pm i}}\circ\psi_{i}^1$, $\phi_{i}^2=\phi_{F_{\pm i}}\circ\psi_{i}^2$.

   On peut choisir, et on choisit, un \'el\'ement $g_{0}\in G_{\tilde{\theta}}(\bar{F})$ et une bijection
$$\iota:\bigsqcup_{i\in I}\Phi_{i}\to \{1,...,d\}\setminus\{(d+1)/2\}$$ 
v\'erifiant les conditions ci-dessous. Pour $i\in I$, on pose $K_{i}=\iota(\Phi_{i})$. Si $i\in I^*$, on pose $k_{i}=\iota(\phi_{i})$. Si $i\in I\setminus I^*$, on pose $K_{i}^{b}=\iota(\Phi_{i}^{b})$ et $k_{i}^{b}=\iota(\phi_{i}^{b})$ pour $b=1,2$. Alors

$\bullet$ $ \bigsqcup_{i\in I^-}K_{i}=\{1,...,(d^--1)/2\}\cup \{d-(d^--3)/2,...d\}$;

$\bullet$ $ \bigsqcup_{i\in I^+}K_{i}=\{d^-+1)/2,...d-(d^--1)/2\}\setminus\{(d+1)/2\}$;

$ \bullet$ pour $i\in I$ et $\phi\in \Phi_{i}$, $\iota(\phi\circ\tau_{i})=d+1-\iota(\phi)$.

$\bullet$ pour $i\in I^*$, $ k_{i}<(d+1)/2$;  

$\bullet$ pour $i\in I\setminus I^*$ et $k\in K_{i}$, on a $k<(d+1)/2$ si $k\in K_{i}^1$ et $k>(d+1)/2$ si $k\in K_{i}^2$;

$\bullet$ $g_{0}T^{\diamond}g_{0}^{-1}=T$;

$\bullet$ pour $t^{\diamond}=(t^{\diamond}_{D},(t^{\diamond}_{i})_{i\in I})\in T(F)$, posons $t=g_{0}t^{\diamond}g_{0}^{-1}$; alors on a $t_{(d+1)/2}=t^{\diamond}_{D}$ et, pour $i\in I$ et $\phi\in \Phi_{i}$, $t_{\iota(\phi)}=\phi(t^{\diamond}_{i})$.

 Notons $\Omega$ le sous-groupe des \'el\'ements du groupe de Weyl de $T$ dans $GL_{d}$ qui sont fix\'es par $\theta$. Il s'identifie au groupe des permutations $\omega$ de l'ensemble $\{1,...,d\}$ telles que $\omega(d+1-k)=d+1-\omega(k)$ pour tout $k$. Posons $\Gamma=Gal(\bar{F}/F)$. Ce groupe agit naturellement sur $T^{\diamond}$. Via l'isomorphisme $t^{\diamond}\mapsto g_{0}t^{\diamond}g_{0}^{-1}$, cette action se transporte en une action sur $T$. Pour simplifier les notations, on note $(\sigma,t)\mapsto \sigma(t)$ cette action. Il existe un unique homomorphisme $\sigma\mapsto \omega(\sigma)$ de $\Gamma$ dans $\Omega$ de sorte que, pour tous $\sigma\in \Gamma$, $t\in T$, et $k\in \{1,...,d\}$, on ait $(\sigma(t))_{k}=\sigma(t_{\omega(\sigma)^{-1}(k)})$. Pour simplifier, on pose simplement $\sigma(k)=\omega(\sigma)(k)$. Evidemment, pour $i\in I$ et $\phi\in \Phi_{i}$, on a $\sigma(\iota(\phi))= \iota(\sigma\circ\phi)$. L'action que l'on vient de d\'efinir  de $\Gamma$ sur $T$ conserve $T_{\theta}$. L'action d\'eduite sur l'ensemble des caract\`eres de $T_{\theta}$ conserve $\Sigma_{\theta}$.

 On d\'efinit une action de $\Gamma$ sur $\hat{T}$,  par $(\sigma(\hat{t}))_{k}=\hat{t}_{\sigma^{-1}(k)}$. Cette action conserve $\hat{T}_{\hat{\theta}}$ et l'action d\'eduite sur l'ensemble des caract\`eres de $\hat{T}_{\hat{\theta}}$ conserve $\check{\Sigma}_{\hat{\theta}}$.

Dans la suite, $T$, $\hat{T}$ et l'ensemble $\{1,...,d\}$ seront toujours munis des actions galoisiennes  que l'on vient de d\'efinir. Pour \'eviter les confusions, on note $T(\bar{F})^{\Gamma}$ plut\^ot que $T(F)$ le sous-groupe des points fixes par $\Gamma$ dans $T(\bar{F})$. On pose $t_{\tilde{x}}=g_{0}t_{\tilde{x}}^{\diamond}g_{0}^{-1}$. Remarquons que $t_{\tilde{x}}$  appartient \`a $T(\bar{F})^{\Gamma}$. 
 
 Pour tout $\beta\in \Sigma_{\theta}$, on note $\Gamma_{\beta}$ le fixateur de $\beta$ dans $\Gamma$, $\Gamma_{\pm \beta}$ le sous-groupe des \'el\'ements de $\Gamma$ qui conservent $\{\pm \beta\}$,  et $F_{\beta}$, resp. $F_{\pm \beta}$, le sous-corps des points fixes par $\Gamma_{\beta}$, resp. $\Gamma_{\pm \beta}$, dans $\bar{F}$. Fixons des $a$-data $(a_{\beta})_{\beta\in \Sigma_{\theta}}$ et des $\chi$-data
$(\chi_{\beta})_{\beta\in \Sigma_{\theta}}$, cf. [LS] paragraphes 2.2 et 2.5.  Pour $\beta\in \Sigma_{\theta}$, $a_{\beta}$ est un \'el\'ement de $F_{\beta}^{\times}$. On a $a_{-\beta}=-a_{\beta}$ et $a_{\sigma(\beta)}=\sigma(a_{\beta})$ pour tout $\sigma\in \Gamma$. Le terme $\chi_{\beta}$ est un caract\`ere de $F_{\beta}^{\times}$. On a  $\chi_{-\beta}=\chi_{\beta}^{-1}$ et $\chi_{\sigma(\beta)}=\chi_{\beta}\circ\sigma^{-1}$ pour tout $\sigma\in \Gamma$. Si $[F_{\beta}:F_{\pm \beta}]=2$, la restriction de $\chi_{\beta}$ \`a $F_{\pm \beta}^{\times}$ est le caract\`ere quadratique $sgn_{F_{\beta}/F_{\pm \beta}}$.

Remarquons que la bijection entre $\Sigma_{\theta}$ et $\check{\Sigma}_{\hat{\theta}}$ est $\Gamma$-\'equivariante. Par cette bijection, on peut, si on en a besoin, remplacer les ensembles d'indices de nos $a$-data et $\chi$-data par $\check{\Sigma}_{\hat{\theta}}$. On peut aussi d\'efinir $\Gamma_{\beta}$,$F_{\beta}$ etc... pour $\beta\in \check{\Sigma}_{\hat{\theta}}$.

\bigskip

\subsection{Le facteur de transfert}

Le facteur $\Delta_{H,\tilde{G}}(y,\tilde{x})$ est un produit de trois termes. Un terme $\Delta_{I}(y,\tilde{x})$ sur lequel nous reviendrons plus tard.  Fixons un sous-ensemble $\underline{\Sigma}$ de repr\'esentants des orbites de l'action de $\Gamma\times\{1,\theta\}$ dans $\Sigma$. Le terme $\Delta_{II}(y,\tilde{x})$  est lui-m\^eme un produit de termes $\Delta_{II,\alpha}(y,\tilde{x})$, pour $\alpha$ dans un sous-ensemble de $ \underline{\Sigma}$.  C'est le sous-ensemble des $\alpha \in \underline{\Sigma}$ qui v\'erifient l'une des conditions suivantes:

(1) $\alpha$ est de type $R_{1}$ et $(\check{\alpha})_{res}$ n'appartient pas \`a $\check{\Sigma}_{H}$;  

(2) $\alpha$ est de type $R_{2}$ et ni $(\check{\alpha})_{res}$, ni $2(\check{\alpha})_{res}$ n'appartiennent \`a $\check{\Sigma}_{H}$; 

(3) $\alpha$ est de type $R_{3}$ et $(\check{\alpha})_{res}$ appartient \`a $\check{\Sigma}_{H}$.

Posons $\alpha=\alpha_{k,l}$. La condition (1) \'equivaut \`a, ou bien $k\in \{1,...,(d^--1)/2\}\cup\{d-(d^--3)/2,...,d\}$ et $l\in \{d^-+1)/2,...,d-(d^--1)/2\}\setminus\{(d+1)/2\}$, ou bien $l\in \{1,...,(d^--1)/2\}\cup\{d-(d^--3)/2,...,d\}$ et $k\in \{d^-+1)/2,...,d-(d^--1)/2\}\setminus\{(d+1)/2\}$. La condition (2) n'est jamais v\'erifi\'ee. En effet, si  $\alpha$ est de type $R_{2}$, soit $k$, soit $l$ est \'egal \`a $(d+1)/2$. Supposons par exemple $l=(d+1)/2$. Si $k\in  \{1,...,(d^--1)/2\}\cup\{d-(d^--3)/2,...,d\}$, $2(\check{\alpha})_{res}$ appartient \`a $\check{\Sigma}_{H^-}$. Si $k\in  \{d^-+1)/2,...,d-(d^--1)/2\}\setminus\{(d+1)/2\}$, $(\check{\alpha})_{res}$ appartient \`a $\check{\Sigma}_{H^+}$. La condition (3) \'equivaut \`a $k=d+1-l$ et $k\in \{1,...,(d^--1)/2\}\cup\{d-(d^--3)/2,...,d\}$.

Pour $\alpha$ v\'erifiant (1), on a
$$\Delta_{II,\alpha}(y,\tilde{x})=\chi_{\alpha_{res}}(\frac{N\alpha(t_{\tilde{x}})-1}{a_{\alpha_{res}}}).$$
Pour $\alpha$ v\'erifiant (3), on a
$$\Delta_{II,\alpha}(y,\tilde{x})=\chi_{\alpha_{res}/2}(\frac{N\alpha(t_{\tilde{x}})+1}{2}).$$

{\bf Remarque.} Dans la d\'efinition de [KS] p.36, il n'y a pas de d\'enominateur $2$. Il nous semble qu'il faut modifier cette d\'efinition comme on vient de l'indiquer.

Le facteur $\Delta_{III}(y,\tilde{x})$ est donn\'e en g\'en\'eral par un produit dans des groupes d'hypercohomologie. En fait, notre choix de $T^{\diamond}$ comme commutant d'un sous-tore de $G_{\tilde{\theta}}$ le simplifie. Le cocycle $v$ de [KS] lemme 4.4.A est trivial: avec les notations de cette r\'ef\'erence, on a $\delta^*=m(\delta)=t_{\tilde{x}}^{\diamond}$, $g=1$, et $u(\sigma)=1$ car $G\times\{1,\theta\}$ est d\'eploy\'e, donc $v(\sigma)=gu(\sigma)\sigma(g)^{-1}=1$.  Alors, d'apr\`es [KS] p.63, le facteur $\Delta_{III}(y,\tilde{x})$ se d\'ecrit de la fa\c{c}on suivante. On d\'efinira plus loin un certain cocycle $a^{\diamond}\in H^1(W_{F},\hat{T})$ (le groupe $W_{F}$ agit sur $\hat{T}$ via l'homomorphisme naturel de $W_{F}$ dans $\Gamma$). Par l'isomorphisme de Langlands entre ce groupe de cohomologie et le groupe des caract\`eres continus de $T(\bar{F})^{\Gamma}$, le cocycle $a^{\diamond}$ d\'efinit un caract\`ere ${\bf a}^{\diamond}$ de $T(\bar{F})^{\Gamma}$. Alors
$$\Delta_{III}(y,\tilde{x})={\bf a}^{\diamond}(t_{\tilde{x}}).$$
D\'ecrivons $a^{\diamond}$. On dit qu'un \'el\'ement $\beta\in \check{\Sigma}_{\hat{\theta}}$ est positif si $\beta$ intervient dans l'action de $\hat{T}_{\hat{\theta}}$ dans le sous-groupe unipotent "triangulaire sup\'erieur" de $\hat{G}_{\hat{\theta}}$. D'autre part, le groupe $\Gamma\times \{\pm 1\}$ agit sur $\check{\Sigma}_{\hat{\theta}}$: l'\'el\'ement $-1$ du groupe $\{\pm 1\}$ agit par $\beta\mapsto -\beta$. Fixons un ensemble de repr\'esentants $ \underline{\check{\Sigma}}_{\hat{\theta}}$ des orbites pour cette action. Pour tout $\beta\in \underline{\check{\Sigma}}_{\hat{\theta}}$, on d\'efinit $r_{\beta}:W_{F}\to \hat{T}$ de la fa\c{c}on suivante. Supposons d'abord que $\beta$ soit asym\'etrique, c'est-\`a-dire que $-\beta$ n'appartienne pas \`a l'orbite de $\beta$ pour $\Gamma$. On fixe un ensemble de repr\'esentants $\{w_{n}\}_{n=1,...,N}$ de $W_{F_{\beta}}\backslash W_{F}$.Pour tout $n$, on pose $\beta_{n}=w_{n}^{-1}\beta$.  Pour tout $(n,w)\in\{1,...,N\}\times W_{F}$, il y a un unique couple $(n',v_{n}(w))\in \{1,...,N\}\times W_{F_{\beta}}$ tel que $w_{n}w=v_{n}(w)w_{n'}$. On pose
$$r_{\beta}(w)=\left(\prod_{n=1,...,N;\beta_{n}>0, w^{-1}\beta_{n}<0}( \check{\beta}_{n})(-1)\right)\left(\prod_{n=1,...,N}\check{\beta}_{n}(\chi_{\beta}(v_{n}(w)))\right).$$
Ici $\check{\beta}_{n}$ est la coracine associ\'ee \`a $\beta_{n}$: c'est un homomorphisme de ${\mathbb C}^{\times}$ dans $\hat{T}_{\hat{\theta}}$. Supposons maintenant que $\beta$ soit sym\'etrique, c'est-\`a-dire que $-\beta$ appartienne \`a l'orbite de $\beta$ pour $\Gamma$. On fixe des \'el\'ements $w_{0}, w_{1},...,w_{N}$ de $W_{F}$ de sorte que $w_{0}^{-1}\beta=-\beta$ et que l'application $n\mapsto \beta_{n}=w_{n}^{-1}\beta$ soit une bijection de $\{1,...,N\}$ sur le sous-ensemble des \'el\'ements positifs de l'orbite de $\beta$. Pour $n\in \{1,...,N\}$, on pose $w_{-n}=w_{0}w_{n}$. Pour tout $(n,w)\in \{1,...,N\}\times W_{F}$, il existe un unique couple $(n',v_{n}(w))\in \{\pm 1,...,\pm N\}\times W_{F_{\beta}}$ tel que $w_{n}w=v_{n}(w)w_{n'}$. On pose
$$r_{\beta}(w)=\prod_{n=1,...,N}\check{\beta}_{n}(\chi_{\beta}(v_{n}(w))).$$ 
On pose ensuite
$$r=\prod_{\beta\in \underline{\check{\Sigma}}_{\hat{\theta}}}r_{\beta}.$$
C'est une application de $W_{T}$ dans $\hat{T}_{\hat{\theta}}$. D'autre part, le groupe $\Omega$ s'identifie au groupe de Weyl de $\hat{T}_{\hat{\theta}}$ dans $\hat{G}_{\hat{\theta}}$. Modulo le choix d'un \'epinglage, on d\'efinit une section de Springer $n:\Omega\to \hat{G}_{\hat{\theta}}$.  Introduisons le $L$-groupe de $T$, qui est le produit semi-direct $^LT=\hat{T }\rtimes W_{F}$. On d\'efinit une application
$$m:{^LT}\to {G}_{\hat{\theta}}\times W_{F}$$
par $m(\hat{t},w)=(\hat{t}r(w)n(\omega(w)),w)$ (on note $\omega(w)$ l'image par $\omega$ de l'image de $w$ dans $\Gamma$). C'est un homomorphisme.  Remarquons que le groupe de Weyl $\Omega_{\hat{H}}$  de $\hat{T}_{\hat{\theta}}$ dans $\hat{H}$ est naturellement un sous-groupe de $\Omega$ et que $\omega$ prend ses valeurs dans ce sous-groupe. D'autre part, tout \'element de $\check{\Sigma}_{H}$ est multiple positif d'un \'el\'ement de $\check{\Sigma}_{\hat{\theta}}$, ce qui permet de d\'efinir $\chi_{\beta}$ pour $\beta\in \check{\Sigma}_{H}$: si $\beta=r\beta'$, avec $r>0$ et $\beta'\in \check{\Sigma}_{\hat{\theta}}$, on pose $\chi_{\beta}=\chi_{\beta'}$. On peut alors remplacer dans les constructions ci-dessus le groupe $\hat{G}_{\hat{\theta}}$ par $\hat{H}$. On affecte d'un indice $H$ les objets relatifs \`a $\hat{H}$. On d\'efinit donc un homomorphisme $m_{H}:{^LT}\to \hat{H}\times W_{F}={^LH}$. On peut consid\'erer que $m$  prend ses valeurs dans $\hat{G}\times W_{F}={^LG}$. L'homomorphisme ${^L\xi}\circ m_{H}$ prend aussi ses valeurs dans $^LG$. On v\'erifie  qu'il existe un cocycle $a^{\diamond}$ de $W_{F}$ dans $\hat{T}$ de sorte que
$$^L\xi\circ m_{H}(\hat{t},w)=a^{\diamond}(w)m(\hat{t},w)$$
pour tout $(\hat{t},w)\in {^LT}$.

\bigskip

\subsection{Le terme $\Delta_{I}(y,\tilde{x})$}

Pour tout $i\in I$, choisissons $X_{i}\in F_{i}^{\times}$ tel que $\tau_{i}(X_{i})=-X_{i}$.  Notons $X$ l'\'el\'ement de $\mathfrak{t}^{\flat}(F)$ qui agit par multiplication par $X_{i}$ sur $F_{i}$ pour tout $i\in I$, et par $0$ sur $D$.  Supposons $X$ r\'egulier dans $\mathfrak{g}_{\theta}(F)$. Introduisons un groupe sp\'ecial orthogonal quasi-d\'eploy\'e $H^{_{'}+}$ d'un espace quadratique de dimension $d^++1$ sur $F$.  Le groupe $H'=H^-\times H^{_{'}+}$ est un groupe endoscopique de $G_{\theta}$. Les classes de conjugaison, ou de conjugaison stable, d'\'el\'ements semi-simples r\'eguliers dans les alg\`ebres de Lie des groupes sp\'eciaux orthogonaux se param\`etrent essentiellement comme en 1.3.  En particulier, on
 peut  introduire des \'el\'ements semi-simples r\'eguliers $Y^-\in \mathfrak{h}^-(F)$ et $Y^+\in \mathfrak{h}^{_{'}+}(F)$ dont les classes de conjugaison stable sont param\'etr\'ees respectivement par $(I^-,(F_{\pm i})_{i\in I^-},(F_{i})_{i\in I^-},(X_{i})_{i\in I^-})$ et $(I^+,(F_{\pm i})_{i\in I^+},(F_{i})_{i\in I^+},(X_{i})_{i\in I^+})$. Posons $Y=(Y^-,Y^+)\in \mathfrak{h}'(F)$.   Supposons de plus $X$  assez proche de $0$.  Les classes de conjugaison stable de $exp(Y)$ et de $exp(X)$ se correspondent. On d\'efinit le facteur de transfert $\Delta_{H',G_{\theta}}(exp(Y),exp(X))$ relatif \`a l'\'epinglage de $G_{\theta}$ que l'on a fix\'e. Comme dans le paragraphe pr\'ec\'edent, c'est le produit de trois termes  que l'on peut calculer en utilisant les m\^emes $a$-data et $\chi$-data que l'on a fix\'ees. Il suffit de comparer les d\'efinitions pour voir que les deux facteurs $\Delta_{I}(exp(Y),exp(X))$ et $\Delta_{I}(y,\tilde{x})$ sont \'egaux. Or on a calcul\'e le premier dans [W] chapitre X pour un certain choix de $a$-data. D\'ecrivons le r\'esultat. Posons $X^T=g_{0}Xg_{0}^{-1}$ (en notant comme une conjugaison l'action adjointe). On  choisit pour $a$-data la famille d\'efinie par $a_{\beta}=\beta(X^T)$ pour tout $\beta\in \Sigma_{\theta}$. On note $Q_{X}$ le polyn\^ome caract\'eristique de $X$ agissant dans $V$. Alors
$$(1)\qquad \Delta_{I}(y,\tilde{x})=\prod_{i\in I^{-*}}sgn_{F_{i}/F_{\pm i}}(\eta c_{i}Q'_{X}(X_{i})),$$
cf. [W] proposition X.8.

{\bf Remarques.}  Dans cette r\'ef\'erence, on avait calcul\'e  le facteur $\Delta_{H',G_{\theta}}(exp(Y),exp(X))$ tout entier, en supposant que $F$ \'etait non-archim\'edien. Mais, avec notre choix de $a$-data, le facteur $\Delta_{II}(exp(Y),exp(X))$ est trivial. Pour $F$ non-archim\'edien, le facteur $\Delta_{III}(exp(Y),exp(X))$ est lui-aussi trivial pour $X$ proche de $0$. La d\'emonstration consistait donc au seul calcul de $\Delta_{I}(exp(Y),exp(X))$. Le calcul de ce terme est le m\^eme, que $F$ soit non-archim\'edien ou r\'eel. Par ailleurs, le terme  que l'on a not\'e ici $\eta$ est le m\^eme que dans [W], bien que sa d\'efinition soit un peu diff\'erente. Le terme que l'on a not\'e $c_{i}$ est \'egal \`a $[F_{i}:F]^{-1}c_{i}$ dans les notations de [W].  

On a suppos\'e $X$ proche de $0$. On peut lever cette hypoth\`ese en remarquant que la formule ci-dessus est insensible au remplacement de $X$ par $\lambda^2X$ pour $\lambda\in F^{\times}$.   Soit $i\in I$. Puisque $y_{i}\tau_{i}(y_{i})=1$, l'\'el\'ement 
$$(2) \qquad X_{i}=(y_{i}-1)(1+y_{i})^{-1}$$
 v\'erifie la propri\'et\'e requise $\tau_{i}(X_{i})=-X_{i}$. On peut  choisir, et on choisit, cet \'el\'ement $X_{i}$ dans les constructions ci-dessus. On v\'erifie que la forte  r\'egularit\'e de $\tilde{x}$ entra\^{\i}ne la r\'egularit\'e de $X$. Remarquons que l'on a la formule d'inversion 
$$(3) \qquad y_{i}=(1+X_{i})(1-X_{i})^{-1}.$$

\bigskip

\subsection{Description de certaines orbites}

 Pour $i\in I^*$, notons $\Gamma_{i}$, resp. $\Gamma_{\pm i}$, le fixateur de $F_{i}$, resp. $F_{\pm i}$ dans $\Gamma$ et \'etendons l'\'el\'ement $\tau_{i}$ de $Gal(F_{i}/F_{\pm i})$ en un \'el\'ement de $\Gamma$. On a donc $\Gamma_{\pm i}=\Gamma_{i}\cup \Gamma_{i}\tau_{i}$.

Fixons $i\in I^-$, $j\in I^+$.  Posons
$$\Sigma(i,j)=\{\alpha_{k,l}; k\in K_{i}, l\in K_{j}\}\cup \{\alpha_{l,k}; k\in K_{i},l\in K_{j}\},$$
$$\Sigma(i,j)_{res}=\{\alpha_{res}; \alpha\in \Sigma(i,j)\}.$$
Les orbites dans $\Sigma(i,j)_{res}$ pour l'action de $\Gamma$ s'identifient aux orbites dans $\Sigma(i,j)$ pour l'action de $\Gamma\times\{1,\theta\}$, ou encore aux orbites dans $K_{i}\times K_{j}$ pour l'action de $\Gamma$. Remarquons que l'on a l'\'egalit\'e  $-(\alpha_{k,l})_{res}=(\alpha_{d+1-k,d+1-l})_{res}$. 

Si $i\not\in I^{-*}$ ou $j\not\in I^{+*}$, $\Sigma(i,j)_{res}$ est form\'e d'orbites asym\'etriques. En effet, si par exemple $i\not\in I^{-*}$, $k$ et $d+1-k$ appartiennent l'un \`a $K_{i}^1$, l'autre \`a $K_{i}^2$, ils ne peuvent donc pas appartenir \`a la m\^eme orbite galoisienne. 

Supposons $i\in I^{-*}$ et $j\in I^{+*}$. Soit $\Xi(i,j)$ un ensemble de repr\'esentants de $\Gamma_{i}\backslash\Gamma/\Gamma_{j}$. On v\'erifie que $((k_{i},\xi(k_{j})))_{\xi\in \Xi(i,j)}$ est une famille de repr\'esentants des orbites de $\Gamma$ dans $K_{i}\times K_{j}$. Soit $\xi\in \Xi(i,j)$. Il y a un \'el\'ement $\xi'\in \Xi(i,j)$ tel que $\tau_{i}\xi\tau_{j}\in \Gamma_{i}\xi'\Gamma_{j}$. Alors la $\Gamma$-orbite de $(\alpha_{k_{i},\xi(k_{j})})_{res}$ est sym\'etrique si $\xi'=\xi$ et asym\'etrique si $\xi'\not=\xi$. Cela r\'esulte des \'egalit\'es
 $$(1) \qquad -(\alpha_{k_{i},\xi(k_{j})})_{res}=(\alpha_{d+1-k_{i},d+1-\xi(k_{j})})_{res}=(\alpha_{\tau_{i}^{-1}(k_{i}),\xi\tau_{j}(k_{j})})_{res}=\tau_{i}^{-1}(\alpha_{k_{i},\tau_{i}\xi\tau_{j}(k_{j})})_{res}.$$
 On fixe un sous-ensemble $\Xi(i,j)_{\pm}$ de $\Xi(i,j)$ tel que, pour tout $\xi\in \Xi(i,j)$, l'intersection  $\Xi(i,j)_{\pm}\cap \{\xi,\xi'\}$ ait exactement un \'el\'ement. Alors l'ensemble des $(\alpha_{k_{i},\xi(k_{j})})_{res}$ pour $\xi\in \Xi(i,j)_{\pm}$ est un ensemble de repr\'esentants des orbites dans $\Sigma(i,j)_{res}$ pour l'action de $\Gamma\times \{\pm 1\}$. 
 
 Fixons maintenant $i\in I^-$. Posons
 $$\Sigma(i)=\{\alpha_{k,d+1-k}; k\in K_{i }\},$$
 $$\Sigma(i)_{res}=\{\alpha_{res}; \alpha\in \Sigma(i)\},\,\,\frac{1}{2}\Sigma(i)_{res}=\{\frac{1}{2}\alpha_{res}; \alpha\in \Sigma(i)\}.$$
 Les orbites dans $\Sigma(i)_{res}$ ou $\frac{1}{2}\Sigma(i)_{res}$ pour l'action de $\Gamma$ s'identifient aux orbites dans $\Sigma(i)$ pour l'action de $\Gamma\times\{1,\theta\}$ ou $\Gamma$ (puisque $\theta$ agit trivialement sur $\Sigma(i)$), ou encore aux orbites dans $K_{i} $ pour l'action de $\Gamma$. Il y en a deux si $i\not\in I^{-*}$, une si $i\in I^{-*}$. Pour $i\not\in I^{-*}$, les deux orbites dans $\Sigma(i)_{res}$ ou $\frac{1}{2}\Sigma(i)_{res}$  sont asym\'etriques.  Si $i\in I^{-*}$, l'unique orbite est sym\'etrique.
 
 On a choisi un ensemble $\underline{\Sigma}$ de repr\'esentants des orbites dans $\Sigma$ pour l'action de $\Gamma\times \{1,\theta\}$. Pour $i\in I^{-}$ et $j\in I^+$, on pose $\underline{\Sigma}(i,j)=\underline{\Sigma}\cap \Sigma(i,j)$. Si $i\in I^{-*}$ et $j\in I^{+*}$, on suppose que $\underline{\Sigma}(i,j)=\{\alpha_{k_{i},\xi(k_{j})}; \xi\in \Xi(i,j)\}$. Pour $i\in I^-$, on pose $\underline{\Sigma}(i)=\underline{\Sigma}\cap \Sigma(i)$. Si $i\in I^{-*}$, on suppose que $\underline{\Sigma}(i)=\{\alpha_{k_{i},d+1-k_{i}}\}$.
 
 On peut dans tout cela remplacer les $\alpha_{k,l}$ par les $\check{\alpha}_{k,l}$, en ajoutant des $\check{}$ dans les d\'efinitions. Du c\^ot\'e des $L$-groupes, c'est le groupe $W_{F}$ qui intervient plut\^ot que $\Gamma$. Il convient alors de consid\'erer que, pour $i\in I^{-*}$ et $j\in I^{+*}$, $\Xi(i,j)$ est un ensemble de repr\'esentants de $W_{F_{i}}\backslash W_{F}/W_{F_{j}}$ et, pour $i\in I^*$, $\tau_{i}$ est un \'el\'ement de $W_{F_{\pm i}}$. On a choisi des ensembles $\underline{\check{\Sigma}}_{\check{\theta}}$, resp. $\underline{\check{\Sigma}}_{H}$, de repr\'esentants des orbites dans $\check{\Sigma}_{\check{\theta}}$, resp. $\check{\Sigma}_{H}$, pour l'action de $\Gamma\times \{\pm 1\}$. Pour $i\in I^-$ et $j\in I^+$, posons $\underline{\check{\Sigma}}(i,j)_{res}=\underline{\check{\Sigma}}_{\check{\theta}}\cap \check{\Sigma}(i,j)_{res}$. Si $i\in I^{-*}$ et $j\in I^{+*}$, on suppose que $\underline{\check{\Sigma}}(i,j)_{res}=\{(\check{\alpha}_{k_{i},\xi(k_{j})})_{res}; \xi\in \Xi(i,j)_{\pm}\}$. Pour $i\in I^-$, posons $\underline{\check{\Sigma}}(i)_{res}=\underline{\check{\Sigma}}_{H}\cap \check{\Sigma}(i)_{res}$. On suppose que $\underline{\check{\Sigma}}_{\check{\theta}}\cap \frac{1}{2}\check{\Sigma}(i)_{res}=\{\frac{1}{2}\beta; \beta\in \underline{\check{\Sigma}}(i)_{res}\}$. Si $i\in I^{-*}$, on suppose $\underline{\check{\Sigma}}(i)_{res}=\{(\alpha_{k_{i},d+1-k_{i}})_{res}\}$.

\bigskip

 \subsection{Choix de $\chi$-data}

  Pour tout $i\in I^{*}$,   fixons un caract\`ere $\chi_{i}$ de $F_{i}^{\times}$ dont la restriction \`a $F_{\pm i}$ co\"{\i}ncide avec $sgn_{F_{i}/F_{\pm i}}$.

Fixons $i\in I^-$, $j\in I^+$.   Si $i\not\in I^{-*}$ ou $j\not\in I^{+*}$,  $\Sigma(i,j)_{res}$ est form\'e d'orbites asym\'etriques. On choisit $\chi_{\beta}=1$ pour tout $\beta\in \Sigma(i,j)_{res}$. Supposons $i\in I^{-*}$ et $j\in I^{+*}$.   Soit $\xi\in \Xi(i,j)$, posons $\beta=(\alpha_{k_{i},\xi(k_{j})})_{res}$. Alors $F_{\beta}$ est le compos\'e des extensions $F_{i}$ et $\xi(F_{j})$. On pose $\chi_{\beta}=\chi_{i}\circ Norm_{F_{\beta}/F_{i}}$. Ensuite, pour $\sigma\in \Gamma$, on pose $\chi_{\sigma(\beta)}=\chi_{\beta}\circ \sigma^{-1}$. Montrons que ces choix conviennent. On doit v\'erifier 

(1) pour $\beta\in \Sigma(i,j)_{res}$, $\chi_{-\beta}=\chi_{\beta}^{-1}$;

(2) soit $\beta\in \Sigma(i,j)_{res}$, supposons l'orbite de $\beta$ sym\'etrique; alors la restriction de $\chi_{\beta}$ \`a $F_{\pm \beta}$ est le caract\`ere $sgn_{F_{\beta}/F_{\pm\beta}}$.

On peut se limiter aux $\beta$ de la forme $\beta=(\alpha_{k_{i},\xi(k_{j})})_{res}$ pour $\xi\in \Xi(i,j)$.   Introduisons $\xi'$ comme dans le paragraphe pr\'ec\'edent et fixons $\gamma\in \Gamma_{i}$ tel que $\tau_{i}\xi\tau_{j}\in \gamma\xi'\Gamma_{j}$. Posons $\beta'=(\alpha_{k_{i},\xi'(k_{j})})_{res}$. L'\'egalit\'e 2.6(1) nous dit que $-\beta=\tau_{i}^{-1}\gamma\beta'$. Donc
$$\chi_{-\beta}=\chi_{\beta'}\circ \gamma^{-1}\tau_{i}=\chi_{i}\circ Norm_{F_{\beta'}/F_{i}}\circ\gamma^{-1}\tau_{i}.$$
 Pour deux extensions $F'\subset F''$ de $F$ et pour $\delta\in \Gamma$, on a l'\'egalit\'e $Norm_{F''/F'}\circ\delta=\delta\circ Norm_{\delta^{-1}(F'')/\delta^{-1}(F')}$. On en d\'eduit ici
 $$\chi_{-\beta}=\chi_{i}\circ\tau_{i}\circ Norm_{F_{\beta}/F_{i}}.$$
 Mais $\chi_{i}\circ\tau_{i}=\chi_{i}^{-1}$ et l'\'egalit\'e pr\'ec\'edente entra\^{\i}ne (1).

 Supposons l'orbite de $\beta$ sym\'etrique. Le corps $F_{\pm\beta}$ contient $F_{\pm i}$: un \'el\'ement de $\Gamma_{\pm \beta}$ envoie forc\'ement $k_{i}$ sur $k_{i}$ ou $d+1-k_{i}$, donc appartient \`a $\Gamma_{\pm i}$.  Le corps $F_{\beta}$ est le compos\'e de $F_{i}$ et de $F_{\pm \beta}$: un \'el\'ement de $\Gamma_{\pm \beta}$ appartient \`a $\Gamma_{\beta}$ si et seulement s'il fixe la premi\`ere composante $k_{i}$ de $\beta$, c'est-\`a-dire si et seulement s'il appartient \`a $\Gamma_{i}$. Puisque l'orbite de $\beta$ est sym\'etrique, $F_{\beta}$ est une extension quadratique de $F_{\pm \beta}$, donc $F_{i}$ et $F_{\pm\beta}$ sont des extensions disjointes de $F_{\pm i}$.  Alors la restriction de $\chi_{\beta}$ \`a $F_{\pm \beta}$ est \'egal \`a $(\chi_{i})_{\vert F_{\pm i}}\circ Norm_{F_{\pm \beta}/F_{\pm i}}$, c'est-\`a-dire \`a $sgn_{F_{i}/F_{\pm i}}\circ Norm_{F_{\pm \beta}/F_{\pm i}}$. Il est connu que c'est bien le caract\`ere $sgn_{F_{\beta}/F_{\pm \beta}}$, d'o\`u (2).

Fixons maintenant seulement $i\in I^{-}$.  Si $i\not\in I^{-*}$, les deux $\Gamma$-orbites dans $\frac{1}{2}\Sigma(i)_{res}$ sont asym\'etriques. On  choisit $\chi_{\beta}=1$ pour tout \'el\'ement $\beta$ de cet ensemble. Si $i\in I^{-*}$, $\frac{1}{2}\Sigma(i)_{res}$ est form\'e d'une unique orbite sym\'etrique. Posons $\beta=\frac{1}{2}(\alpha_{k_{i},d+1-k_{i}})_{res}$. Alors $F_{\beta}=F_{i}$, $F_{\pm \beta}=F_{\pm i}$. On choisit $\chi_{\beta}=\chi_{i}$ et $\chi_{\sigma(\beta)}=\chi_{i}\circ\sigma^{-1}$ pour tout $\sigma\in \Gamma$.

\bigskip

\subsection{Calcul de $\Delta_{II}(y,\tilde{x})$}

Consid\'erons l'ensemble  des  \'el\'ements de $ \Sigma$ qui v\'erifient la condition 2.4(1). C'est la r\'eunion sur les couples $(i,j)\in I^-\times I^+$ des ensembles $\Sigma(i,j)$. Nos choix de $\chi$-data et d'ensemble $\underline{\Sigma}$ entra\^{\i}nent que la contribution \`a $\Delta_{II}(y,\tilde{x})$ des \'el\'ements de $\Sigma$ v\'erifiant la condition 2.4(1) est
$$ \prod_{i\in I^{-*},j\in I^{+*}}\prod_{\xi\in \Xi(i,j)}\Delta_{II,\alpha_{k_{i},\xi(k_{j})}}(y,\tilde{x}).$$
Fixons $i,j,\xi$ intervenant dans ce produit. Posons $\beta=(\alpha_{k_{i},\xi(k_{j})})_{res}$. On a choisi 
$$a_{\beta}=\beta(X^T)=X^T_{k_{i}}-X^T_{\xi(k_{j})}=X_{i}-\xi(X_{j}).$$
 On a 
$$\alpha_{k_{i},\xi(k_{j})}(t_{\tilde{x}})=t_{\tilde{x},k_{i}}t_{\tilde{x},\xi(k_{j})}^{-1}=c_{i}\tau_{i}(x_{i})^{-1}\xi(c_{j}^{-1}\tau_{j}(x_{j})),$$
$$(\theta(\alpha_{k_{i},\xi(k_{j})}))(t_{\tilde{x}})=\alpha_{\xi\tau_{j}(k_{j}),\tau_{i}(k_{i})}(t_{\tilde{x}})=\xi(c_{j}x_{j}^{-1})c_{i}^{-1}x_{i}.$$
On a $x_{i}/\tau_{i}(x_{i})=y_{i}$ et $x_{j}/\tau_{j}(x_{j})=y_{j}$. Donc
 $$N\alpha_{k_{i},\xi(k_{j})}(t_{\tilde{x}})=y_{i }\xi(y_{j})^{-1}.$$
On a $\chi_{\beta}=\chi_{i}\circ Norm_{F_{\beta}/F_{i}}$. D'o\`u
$$\Delta_{II,\alpha_{k_{i},\xi(k_{j})}}(y,\tilde{x})=\chi_{i}\circ Norm_{F_{\beta}/F_{i}}((y_{i}\xi(y_{j})^{-1}-1)(X_{i}-\xi(X_{j}))^{-1}).$$
Pour $\lambda\in F_{\beta}$, on a 
$$Norm_{F_{\beta}/F_{i}}(\lambda)=\prod_{\sigma\in \Gamma_{i}/(\Gamma_{i}\cap \xi\Gamma_{j}\xi^{-1})}\sigma(\lambda).$$
Puisque $y_{i}$ et $X_{i}$ sont fixes par $\Gamma_{i}$, on obtient
$$\Delta_{II,\alpha_{k_{i},\xi(k_{j})}}(y,\tilde{x})=\chi_{i}(\prod_{\sigma\in \Gamma_{i}/(\Gamma_{i}\cap \xi\Gamma_{j}\xi^{-1})}\sigma\xi(y_{j})^{-1}(y_{i}-\sigma\xi(y_{j}))(X_{i}-\sigma\xi(X_{j}))^{-1}).$$
On doit faire le produit de ces expressions quand $\xi$ d\'ecrit $\Xi(i,j)$. Mais quand $\xi$ d\'ecrit $\Xi(i,j)$ et $\sigma$ d\'ecrit $\Gamma_{i}/(\Gamma_{i}\cap \xi\Gamma_{j}\xi^{-1})$, le produit $\sigma\xi$ d\'ecrit $\Gamma/\Gamma_{j}$. Remarquons que le nombre d'\'el\'ements de cet ensemble est $[F_{j}:F]$ et que, puisque $Norm_{F_{j}/F_{\pm j}}(y_{j})=1$, a fortiori $Norm_{F_{j}/F}(y_{j})=1$. On obtient
$$\prod_{\xi\in \Xi(i,j)}\Delta_{II,\alpha_{k_{i},\xi(k_{j})}}(y,\tilde{x})=\chi_{i}( P_{j}(y_{i})Q_{j}(X_{i})^{-1}),$$
o\`u $P_{j}$, resp. $Q_{j}$, est le polyn\^ome caract\'eristique de $y_{j}$, resp. $X_{j}$, agissant sur $F_{j}$.  La contribution \`a $\Delta_{II}(y,\tilde{x})$ des \'el\'ements de $\Sigma$ v\'erifiant la condition 2.4(1) est donc
$$(1) \qquad \prod_{i\in I^{-*},j\in I^{+*}}\chi_{i}( P_{j}(y_{i})Q_{j}(X_{i})^{-1}).$$

Consid\'erons maintenant l'ensemble des \'el\'ements de $\Sigma$ qui v\'erifient la condition 2.4(3). C'est la r\'eunion sur les $i\in I^{-}$ des ensembles $ \Sigma(i)$. Nos choix de $\chi$-data et d'ensemble $\underline{\Sigma}$ entra\^{\i}nent que la contribution \`a $\Delta_{II}(y,\tilde{x})$ des \'el\'ements de $\Sigma$ v\'erifiant la condition 2.4(3) est
$$ \prod_{i\in I^{-*}}\Delta_{II,\alpha_{k_{i},d+1-k_{i}}}(y,\tilde{x}).$$
  Soit $i\in I^{-*}$. On a $N\alpha_{k_{i},d+1-k_{i}}=\alpha_{k_{i},d+1-k_{i}}$, et on calcule comme ci-dessus
$N\alpha_{k_{i},d+1-k_{i}}(t_{\tilde{x}})=y_{i}$. Posons $\beta=\frac{1}{2}(\alpha_{k_{i},d+1-k_{i}})_{res}$. On a choisi $\chi_{\beta}=\chi_{i}$ et on obtient
$$\Delta_{II,\alpha_{k_{i},d+1-k_{i}}}(y,\tilde{x})=\chi_{i}((y_{i}+1)/2).$$
La contribution \`a $\Delta_{II}(y,\tilde{x})$ des \'el\'ements de $\Sigma$ qui v\'erifient la condition 2.4(3) est donc
$$(2) \qquad \prod_{i\in I^{-*}}\chi_{i}((y_{i}+1)/2).$$

\bigskip

\subsection{D\'ebut du calcul de $\Delta_{III}(y,\tilde{x})$}

Parce que $T$ est un produit de tores induits, l'isomorphisme  entre   $H^1(W_{F},\hat{T})$ et  le groupe des caract\`eres de $T(\bar{F})^{\Gamma}$ s'explicite ais\'ement. D\'ecrivons le r\'esultat dans le cas qui nous int\'eresse. Pour $i\in I^*$, l'application
$$\begin{array}{ccc}W_{F_{i}}&\to &{\mathbb C}^{\times}\\ w&\mapsto &a^{\diamond}(w)_{k_{i}}\\ \end{array}$$
est un caract\`ere de $W_{F_{i}}$, auquel correspond un caract\`ere ${\bf a}_{i}^{\diamond}$ de $F_{i}^{\times}$. Pour $i\in I\setminus I^*$, on a deux applications
$$\begin{array}{ccc}W_{F_{\pm i}}&\to &{\mathbb C}^{\times}\\ w&\mapsto &a^{\diamond}(w)_{k_{i}^{b}}\\ \end{array}$$
pour $b=1,2$. Ce sont des caract\`eres de $F_{\pm i}^{\times}$, dont se d\'eduisent deux caract\`eres ${\bf a}_{i}^{b,\diamond}$ de $F_{\pm i}^{\times}$. Enfin l'application
$$\begin{array}{ccc}W_{F}&\to &{\mathbb C}^{\times}\\ w&\mapsto &a^{\diamond}(w)_{(d+1)/2}\\ \end{array}$$
est un caract\`ere de $W_{F}$, dont se d\'eduit un caract\`ere ${\bf a}_{D}^{\diamond}$ de $F^{\times}$. On a alors l'\'egalit\'e
$${\bf a}^{\diamond}(t_{\tilde{x}})={\bf a}^{\diamond}_{D}(t_{\tilde{x},(d+1)/2})\left(\prod_{i\in I\setminus I^*}{\bf a}^{1,\diamond}_{i}(t_{\tilde{x},k^1_{i}}){\bf a}_{i}^{2,\diamond}(t_{\tilde{x},k_{i}^2})\right)\left(\prod_{i\in I^*}{\bf a}^{\diamond}_{i}(t_{\tilde{x},k_{i}})\right),$$
c'est-\`a-dire
$$(1) \qquad {\bf a}^{\diamond}(t_{\tilde{x}})={\bf a}^{\diamond}_{D}( c_{D}x_{D}^{-1})\left(\prod_{i\in I\setminus I^*}{\bf a}^{1,\diamond}_{i}( \psi_{i}^1(c_{i}\tau_{i}(x_{i})^{-1})){\bf a}_{i}^{2,\diamond}( \psi_{i}^2(c_{i}\tau_{i}(x_{i})^{-1}))\right)$$
$$\left(\prod_{i\in I^*}{\bf a}^{\diamond}_{i}( c_{i}\tau_{i}(x_{i})^{-1})\right).$$

  Explicitons le cocycle $a^{\diamond}$. Le lemme X.4 de [W] calcule la section de Springer $n$ relative \`a un \'epinglage convenable. Soit $w\in W_{F}$. Identifions  $\omega(w)$ \`a une matrice de permutation. Pour $k\in \{1,...,d\}$, notons $s(k,w)$ le nombre d'\'el\'ements $l\in \{1,...,d\}$ tels que $l<k$ et $w^{-1}(l)>w^{-1}(k)$. Notons $S(w)$ la matrice diagonale de coefficients $S(w)_{k}=(-1)^{s(k,w)}$. Alors $n(\omega(w))=S(w)\omega(w)$. On a une formule similaire pour $n_{H}(\omega(w))$. Quand on tient compte de la fa\c{c}on dont on a plong\'e $\hat{H}$ dans $\hat{G}$, on obtient le r\'esultat suivant. Pour $k\in \{1,...,d\}$, notons $s_{H}(k,w)$ le nombre d'\'el\'ements de l'ensemble suivant:
 
 - si $k\in \{1,..., (d^--1)/2\}\cup\{d-(d^--3)/2,...,d\}$, c'est l'ensemble des $l\in \{1,...,(d^--1)/2\}\cup\{d-(d^--3)/2,...,d\}$ tels que $l<k$ et $w^{-1}(l)>w^{-1}(k)$;
 
 - si $k\in \{(d^-+1)/2,...,d-(d^--1)/2\}$, c'est l'ensemble des $l\in \{(d^-+1)/2,...,d-(d^--1)/2\}$ tels que $l<k$ et $w^{-1}(l)>w^{-1}(k)$.
 
 Notons $S_{H}(w)$ la matrice diagonale de coefficients $S_{H}(w)_{k}=(-1)^{s_{H}(k,w)}$. Alors $^L\xi( n_{H}(\omega(w)))=S_{H}(w)\omega(w)$. La formule de 2.4 qui d\'efinit le cocycle $a^{\diamond}$ se r\'ecrit donc
 $$^L\xi(r_{H}(w))S_{H}(w)\omega(w)\rho(w)=a^{\diamond}(w)r(w)S(w)\omega(w),$$
 ou encore, puisque $\rho(w)$ et $\omega(w)$ commutent,
 $$ a^{\diamond}(w)={^L\xi}(r_{H}(w))r(w)^{-1}S_{H}(w)S(w)^{-1}\rho(w).$$
 Les calculs ci-dessus permettent d\'ej\`a d'expliciter la matrice $S_{H}(w)S(w)^{-1}$.  Notons $e(w )$ le nombre d'\'el\'ements de l'ensemble des $l\in \{1,...,(d^--1)/2\}$ tels que $w^{-1}(l)\geq d-(d^--3)/2$. Alors   $S_{H}(w)S(w)^{-1}$ est la matrice diagonale de coefficients
 $$(2) \qquad (S_{H}(w)S(w)^{-1})_{k}=\left\lbrace\begin{array}{cc}1&\text{ si }k\leq (d^--1)/2,\\ (-1)^{e(w)},&\text{ si }(d^-+1)/2\leq k\leq d-(d^--1)/2,\\ 1,&\text{ si }k\geq d-(d^--3)/2\text{ et }w^{-1}(k)\geq d-(d^--3)/2,\\ -1,&\text{ si }k\geq d-(d^--3)/2\text{ et }w^{-1}(k)\leq (d^--1)/2.\\ \end{array}\right.$$
 
 On a 
 $${^L\xi}(r_{H}(w))r(w)^{-1}=\left(\prod_{\beta\in \underline{\check{\Sigma}}_{H}}{^L\xi}(r_{H,\beta}(w))\right)\left(\prod_{\beta\in \underline{\check{\Sigma}}_{\hat{\theta}}}r_{\beta}(w)^{-1}\right).$$  
 Dans l'intersection $\check{\Sigma}_{H}\cap \check{\Sigma}_{\hat{\theta}}$, on peut supposer que l'on choisit les m\^emes repr\'esentants d'orbites dans les deux produits. On voit alors que les contributions de cette intersection \`a chacun des produits  sont inverses l'une de l'autre. Elles disparaissent. Il reste la contribution des \'el\'ements de $\check{\Sigma}_{\hat{\theta}}$ qui n'appartiennent pas \`a $\check{\Sigma}_{H}$ et celle des \'el\'ements de $\check{\Sigma}_{H}$ qui n'appartiennent pas \`a $\check{\Sigma}_{\hat{\theta}}$. Le second ensemble est la r\'eunion sur les $i\in I^-$ des $ \check{\Sigma}(i)_{res}$. Le premier ensemble est lui-m\^eme r\'eunion de deux ensembles:
 
 - la r\'eunion sur les $(i,j)\in I^-\times I^+$ des $\check{\Sigma}(i,j)_{res}$;
  
 - la r\'eunion sur les $i\in I^-$ des $\frac{1}{2}\check{\Sigma}(i)_{res}$.

 D'apr\`es nos choix d'ensembles de repr\'esentants d'orbites, on obtient
 $${^L\xi}(r_{H}(w))r(w)^{-1}=\left(\prod_{i\in I^-, j\in I^+}\prod_{\beta\in \underline{\check{\Sigma}}(i,j)_{res}}r_{\beta}(w)^{-1}\right)\left(\prod_{i\in I^-}\prod_{\beta\in \underline{\check{\Sigma}}(i)_{res}}{^L\xi}(r_{H,\beta}(w))r_{\beta/2}(w)^{-1}\right).$$
  Le second produit se simplifie. En effet, pour $i\in I^-$ et $\beta\in \underline{\check{\Sigma}}(i)_{res}$,  la seule diff\'erence dans les d\'efinitions de $^L\xi(r_{H,\beta}(w))$ et $r_{\beta/2}(w)$ provient des coracines qui interviennent dans les d\'efinitions: celles intervenant dans le second terme sont les doubles de celles intervenant dans le premier. Donc $r_{\beta/2}(w)={^L\xi}(r_{H,\beta}(w))^2$. D'autre part, soient $i\in I^-\setminus I^{-*}$, $j\in I^+$ et $\beta\in \underline{\check{\Sigma}}(i,j)_{res}$. L'orbite de $\beta$ pour l'action de $\Gamma\times \{\pm 1\}$ est asym\'etrique et l'action de $\Gamma$ pr\'eserve la positivit\'e: un \'el\'ement $(\alpha_{k,l})_{res}$ de l'orbite est positif si $k\in K_{i}^1$ et n\'egatif si $k\in K_{i}^2$. Donc $r_{\beta}(w)=1$ pour tout $w$. De m\^eme, pour $i\in I^{-}\setminus I^{-*}$ et $\beta\in \underline{\check{\Sigma}}(i)_{res}$, $r_{H,\beta}(w)=1$ pour tout $w$. On obtient  
 $$(3) \qquad {^L\xi}(r_{H}(w))r(w)^{-1}=\left(\prod_{i\in I^{-*}, j\in I^+}\prod_{\beta\in \underline{\check{\Sigma}}(i,j)_{res}}r_{\beta}(w)^{-1}\right)\left(\prod_{i\in I^{-*}}\prod_{\beta\in \underline{\check{\Sigma}}(i)_{res}}{^L\xi}(r_{H,\beta}(w))^{-1}\right).$$
 
 \bigskip
 
  \subsection{Le caract\`ere ${\bf a}^{\diamond}_{D}$}
 
 On doit calculer  de $a^{\diamond}(w)_{(d+1)/2}$ pour $w\in W_{F}$. Le terme $\rho(w)_{(d+1)/2}$ vaut $\chi(w)$. Le terme $(S_{H}(w)S(w)^{-1})_{(d+1)/2}$ vaut $(-1)^{e(w)}$. Les termes $r_{\beta}(w)$ et $^L\xi(r_{H,\beta}(w))$ ont une $(d+1)/2$-i\`eme composante \'egale \`a $1$. On obtient  
 $$a^{\diamond}(w)_{(d+1)/2}=\chi(w)(-1)^{e(w)}.$$
 Identifions le caract\`ere $w\mapsto (-1)^{e(w)}$. Pour $i\in I^-$, notons $e_{i}(w)$ le nombre des $k\in K_{i}$ tels que $k\leq (d^--1)/2$ et $w^{-1}(k)\geq d-(d^--3)/2$. Pour $i\not\in I^{-*}$, l'action galoisienne sur $K_{i}$ conserve  la positivit\'e et $e_{i}(w)=0$. Donc
 $$e(w)=\sum_{i\in I^{-*}}e_{i}(w).$$
 Soit $i\in I^{-*}$. On peut choisir des repr\'esentants $(w_{n})_{n=1,...,N}$ de $W_{F_{\pm i}}\backslash W_{F}$ de sorte que l'application $n\mapsto w_{n}^{-1}(k_{i})$ soit une bijection de $\{1,...,N\}$ sur l'ensemble des \'el\'ements $k\in K_{i}$ tels que $k\leq (d^--1)/2$. Pour $n=1,...,N$, soit $(n',v_{n}(w))\in \{1,...,N\}\times W_{F_{\pm i}}$ tel que $w_{n}w=v_{n}(w)w_{n'}$. On a $w^{-1}w_{n}^{-1}(k_{i})\geq d-(d^--3)/2$ si et seulement si $v_{n}(w)\in W_{F_{\pm i}}\setminus W_{F_{i}}$, ou encore si et seulement si $sgn_{F_{i}/F_{\pm i}}(v_{n}(w))=-1$. Donc $(-1)^{e_{i}(w)}=\prod_{n=1,...,N}sgn_{F_{i}/F_{\pm i}}(v_{n}(w))$. Le membre de droite de cette \'egalit\'e d\'efinit le transfert du caract\`ere $sgn_{F_{i}/F_{\pm i}}$ de $W_{F_{\pm i}}$ en un caract\`ere de $W_{F}$. En termes de caract\`eres de $F_{\pm i}^{\times}$ et $F^{\times}$, le transfert se traduit par la restriction. Le caract\`ere $w\mapsto (-1)^{e_{i}(w)}$ de $W_{F}$ correspond donc \`a la restriction \`a $F^{\times}$ du caract\`ere $sgn_{F_{i}/F_{\pm i}}$ de $F_{\pm i}^{\times}$. D'o\`u
 $${\bf a}^{\diamond}_{D}=\chi\prod_{i\in I^{-*}}(sgn_{F_{i}/F_{\pm i}})_{\vert F^{\times}}.$$
 
 \bigskip

 \subsection{Calcul de ${\bf a}^{b,\diamond}_{i}$, pour $i\in I^-\setminus I^{-*}$}
 
   Soit  $i\in I^-\setminus I^{-*}$. Pour $b=1,2$, on doit calculer $a^{\diamond}(w)_{k_{i}^b}$ pour $w\in W_{F_{\pm i}}$. Le terme $\rho(w)_{k_{i}^b}$ vaut $1$. Le terme $(S_{H}(w)S(w)^{-1})_{k_{i}^b}$ aussi:   cela r\'esulte de la formule 2.9(2) et, dans le cas $b=2$, auquel cas $k_{i}^2\geq d-(d^--3)/2$, de l'hypoth\`ese $w\in W_{F_{\pm i}}$ qui entra\^{\i}ne $w^{-1}(k_{i}^2)=k_{i}^2$.  Le $k_{i}^b$-i\`eme coefficient du membre de droite de 2.9(3) vaut aussi $1$. On en d\'eduit
 
   pour $i\in I^-\setminus I^{-*}$, ${\bf a}^{b,\diamond}_{i}=1$.

 \bigskip
 
 \subsection{Calcul de ${\bf a}^{\diamond}_{i}$, pour $i\in  I^{-*}$}
 
 Soit $i\in I^{-*}$. On doit calculer $a^{\diamond}(w)_{k_{i}}$ pour $w\in W_{F_{i}}$. On a $\rho(w)_{k_{i}}=1$. Le terme $(S_{H}(w)S(w)^{-1})_{k_{i}}$ vaut $1$  par le m\^eme argument que dans le paragraphe pr\'ec\'edent. Dans le produit 2.9(3), seuls contribuent au $k_{i}$-i\`eme coefficient ceux pour lequel l'indice $i$ est notre $i$. 
 
 Soient $j\in I^+\setminus I^{+*}$ et $\beta\in \underline{\check{\Sigma}}(i,j)_{res}$. L'orbite de $\beta$ est asym\'etrique, donc $\chi_{\beta}=1$. Avec les notations de 2.4, on a simplement
 $$r_{\beta}(w)^{-1}=\prod_{n=1,...,N, \beta_{n}>0, w^{-1}(\beta_{n})<0}\check{\beta}_{n}(-1).$$
 Ne contribuent \`a la $k_{i}$-i\`eme composante que les $\beta_{n}$ de la forme $(\check{\alpha}_{k_{i},l})_{res}$ ou $(\check{\alpha}_{d+1-k_{i},l})_{res}$.  Les racines de la premi\`ere forme sont positives, celles de la seconde forme sont n\'egatives. Un \'el\'ement $w\in W_{F_{i}}$ fixe $k_{i}$ et $d+1-k_{i}$ donc respecte la positivit\'e. La $k_{i}$-i\`eme composante de $r_{\beta}(w)$ vaut donc $1$ pour $w\in W_{F_{i}}$. 
 
 Soit maintenant $j\in I^{+*}$. On a $ \underline{\check{\Sigma}}(i,j)_{res}=\{(\check{\alpha}_{k_{i},\xi(k_{j})})_{res}; \xi\in \Xi(i,j)_{\pm}\}$. Fixons $\xi\in \Xi(i,j)_{\pm}$, posons $\beta=(\check{\alpha}_{k_{i},\xi(k_{j})})_{res}$. 
 
 Supposons d'abord l'orbite de $\beta$ asym\'etrique. Il intervient dans $r_{\beta}(w)^{-1}$ un produit de $\check{\beta}_{n}(-1)$. Sa $k_{i}$-i\`eme composante vaut $1$ pour la m\^eme raison que ci-dessus. Il reste le produit
 $$\prod_{n=1,...,N}\check{\beta}_{n}(\chi_{\beta}(v_{n}(w)))^{-1}.$$
 Pour $n=1,...,N$, la $k_{i}$-i\`eme composante de $\check{\beta}_{n}(\chi_{\beta}(v_{n}(w)))^{-1}$ vaut
 
 (1)  $\chi_{\beta}(v_{n}(w))^{-1}$ si $\beta_{n}$ est de la forme $(\check{\alpha}_{k_{i},l})_{res}$, c'est-\`a-dire si $w_{n}\in W_{F_{i}}$;
 
 (2) $\chi_{\beta}(v_{n}(w))$ si $\beta_{n}$ est de la forme $(\check{\alpha}_{d+1-k_{i},l})_{res}$, c'est-\`a-dire si $w_{n}\in W_{F_{i}}\tau_{i}$;

  $1$ sinon.

  L'ensemble  $\{w_{n}; n=1,...,N, w_{n}\in W_{F_{i}}\}$ est un ensemble de repr\'esentants de $W_{F_{\beta}}\backslash W_{F_{i}}$. L'application
  $$\begin{array}{ccc}W_{F_{i}}&\to &{\mathbb C}^{\times}\\ w&\mapsto&\prod_{n=1,...,N,w_{n}\in \Gamma_{i}}\chi_{\beta}(v_{n}(w))^{-1}\\ \end{array}$$
  est le transfert de $W_{F_{\beta}}$ \`a $W_{F_{i}}$ du caract\`ere $\chi_{\beta}^{-1}$.  Le transfert correspond \`a la restriction de $F_{\beta}^{\times}$ \`a $F_{i}^{\times}$. Puisque $\chi_{\beta}^{-1}=\chi_{i}^{-1}\circ Norm_{F_{\beta}/F_{i}}$, cette restriction est \'egale \`a $\chi_{i}^{[-F_{\beta}:F_{i}]}$. Le produit des termes (1) est donc $\chi_{i}(w)^{-[F_{\beta}:F_{i}]}$. L'ensemble $\{w_{n}\tau_{i}^{-1}; n=1,...,N,, w_{n}\in W_{F_{i}}\tau_{i}\}$ est encore un ensemble de repr\'esentants de $W_{F_{\beta}}\backslash W_{F_{i}}$. Remarquons que, pour $w_{n}$ dans cet ensemble, l'\'egalit\'e
  $w_{n}w=v_{n}(w)w_{n'}$ \'equivaut \`a $w_{n}\tau_{i}^{-1}(\tau_{i}w\tau_{i}^{-1})=v_{n}(w)w_{n'}\tau_{i}^{-1}$. Le m\^eme raisonnement que ci-dessus montre que le produit des termes (2) vaut $\chi_{i}(\tau_{i}w\tau_{i}^{-1})^{[F_{\beta}:F_{i}]}$. On a l'\'egalit\'e $\chi_{i}\circ\tau_{i}=\chi_{i}^{-1}$ de caract\`eres de $F_{i}^{\times}$, et cela se traduit du c\^ot\'e galoisien par l'\'egalit\'e $\chi_{i}(\tau_{i}w\tau_{i}^{-1})=\chi_{i}(w)^{-1}$. Le produit des termes (2) est donc \'egal \`a celui des termes (1). La contribution totale est donc $\chi_{i}(w)^{-2[F_{\beta}:F_{i}]}$. Ce calcul montre aussi que $[F_{\beta}:F_{i}]$ est \'egal \`a la fois au nombre d'\'el\'ements de la $\Gamma$-orbite de $\beta$ de la forme $(\check{\alpha}_{k_{i},l})_{res}$ et au nombre d'\'el\'ements de la m\^eme orbite de la forme $(\check{\alpha}_{d+1-k_{i},l})_{res}$. Ce deuxi\`eme nombre est \'egal \`a celui des \'el\'ements de la $\Gamma$-orbite de $-\beta$ de la forme $(\check{\alpha}_{k_{i},l})_{res}$. Donc $2[F_{\beta}:F_{i}]$ est le nombre d'\'el\'ements de la $\Gamma\times\{\pm 1\}$-orbite de $\beta$ de la forme $(\check{\alpha}_{k_{i},l})_{res}$.
  
  Supposons maintenant l'orbite de $\beta$ sym\'etrique. Cette fois, les $\beta_{n}$ sont positifs, donc ne peuvent pas \^etre de la forme $(\check{\alpha}_{d+1-k_{i},l})_{res}$. La $k_{i}$-i\`eme composante de $r_{\beta}(w)^{-1}$ est donc le produit des $\chi_{\beta}(v_{n}(w))^{-1}$ sur les $n=1,...,N$ tels que $\beta_{n}$ soit de la forme $(\check{\alpha}_{k_{i},l})_{res}$. Cette condition \'equivaut \`a $w_{n}\in W_{F_{i}}$. Parce que $W_{F_{i}}\cap W_{F_{\pm \beta}}=W_{F_{\beta}}$, on v\'erifie que l'ensemble des $w_{n}$ qui satisfont cette condition est un ensemble de repr\'esentants de $W_{F_{\beta}}\backslash W_{F_{i}}$ et le m\^eme calcul que ci-dessus montre que la contribution de l'orbite de $\beta$ est $\chi_{i}(w)^{-[F_{\beta}:F_{i}]}$. De plus, $[F_{\beta}:F_{i}]$ est le nombre d'\'el\'ements de la forme $(\check{\alpha}_{k_{i},l})_{res}$ dans la $\Gamma$-orbite de $\beta$, ou dans la $\Gamma\times\{\pm 1\}$-orbite, cela revient au m\^eme.
  
  On  a calcul\'e la contribution de chaque $\Gamma\times \{\pm 1\}$ orbite dans $\check{\Sigma}(i,j)$. Faisons le produit de ces contributions. On obtient $\chi_{i}(w)^{-1}$ \'elev\'e \`a la puissance le nombre d'\'el\'ements de $\check{\Sigma}(i,j)_{res}$ de la forme $(\check{\alpha}_{k_{i},l})_{res}$. Ce nombre est \'evidemment \'egal au nombre d'\'el\'ements de $K_{j}$ c'est-\`a-dire $[F_{j}:F]$. La contribution de  $\check{\Sigma}(i,j)$ est donc $\chi_{i}(w)^{-[F_{j}:F]}$.
  
  Il reste encore \`a calculer la contribution  du deuxi\`eme produit de 2.9(3). On a choisi $\underline{\check{\Sigma}}(i)_{res}=\{(\alpha_{k_{i},d+1-k_{i}})_{res}\}$. On obtient imm\'ediatement pour contribution $\chi_{i}(w)^{-1}$. Finalement,
  
  pour $i\in I^{-*}$, ${\bf a}_{i}^{\diamond}=\chi_{i}^{-1-\sum_{j\in I^{+*}}[F_{j}:F]}$.
  
  \bigskip
  
  \subsection{Calcul de ${\bf a}_{j}^{b,\diamond}$ pour $j\in I^{+}\setminus I^{+*}$}
  
 Soit $j\in I^+\setminus I^{+*}$. Pour $b=1,2$, on doit calculer $a^{\diamond}(w)_{k_{j}^b}$ pour $w\in W_{F_{\pm j}}$. On a $\rho(w)_{k_{j}^b}=\chi(w)$.  C'est en fait la restriction de $\chi$ \`a $W_{F_{\pm j}}$ qui intervient ici. En termes de caract\`ere de $F_{\pm j}^{\times}$, il s'agit de $\chi\circ Norm_{F_{\pm j}/F}$. On a $(S_{H}(w)S(w)^{-1})_{k_{j}^b}=(-1)^{e(w)}$. 
 
 Dans le produit 2.9(3), seuls contribuent au $k_{j}^b$-i\`eme coefficient ceux pour lequel l'indice $j$ est notre $j$. Fixons $i\in I^{-*}$.  Puisque $j\not\in I^{+*}$, toutes les orbites dans $\check{\Sigma}(i,j)_{res}$ sont asym\'etriques. Notons $E$ la r\'eunion des $\Gamma$-orbites  des $\beta\in \underline{\check{\Sigma}}(i,j)_{res}$. C'est un sous-ensemble de $\check{\Sigma}(i,j)$, stable par $\Gamma$ et tel que $\check{\Sigma}(i,j)=E\sqcup(-E)$. La contribution du couple $(i,j)$ au produit 2.9(3) est
 $$\prod_{\beta\in E; \beta>0, w^{-1}(\beta)<0}\check{\beta}(-1).$$
 Notons $E'$, resp. $E''$, le sous-ensemble des \'el\'ements de $E$ de la forme $(\check{\alpha}_{k,k_{j}^b})_{res}$, resp. $(\check{\alpha}_{k,d+1-k_{j}^b})_{res}$. Notons $e'(w)$, resp. $e''(w)$, le nombre d'\'el\'ements $\beta\in E'$, resp. $\beta\in E''$, tels que $\beta>0$ et $w^{-1}(\beta)<0$.  La $k_{j}^b$-i\`eme coordonn\'ee du produit ci-dessus est $(-1)^{e'(w)+e''(w)}$. Notons $K'$, resp. $K''$, l'ensemble des $k\in K_{i}$ tels que  $(\check{\alpha}_{k,k_{j}^b})_{res}\in E'$, resp.  $(\check{\alpha}_{k,d+1-k_{j}^b})_{res}\in E''$. Ce sont des sous-ensembles invariants par $\Gamma_{F_{\pm j}}$.  Posons $K_{+}=K_{i}\cap \{1,...,(d^--1)/2\}$ et $K_{-}=K_{i}\cap \{d-(d^--3)/2,...,d\}$.Si $\beta=(\check{\alpha}_{k,l})_{res}$, avec $k\in K_{i}$ et $l\in K_{j}$, la condition $\beta>0$, resp. $\beta<0$, \'equivaut \`a $k \in K_{+}$, resp. $k\in K_{-}$. Donc $e''(w)$ est le nombre d'\'el\'ements de l'ensemble des $k\in K''\cap K_{+}$ tels que $w^{-1}(k)\in K''\cap K_{-}$.  Pour deux signes $\epsilon,\epsilon'=\pm$, notons $K''_{\epsilon,\epsilon'}(w)$ l'ensemble des $k\in K''\cap K_{\epsilon}$ tels que $w^{-1}(k)\in K''\cap K_{\epsilon'}$. On a
 $$w^{-1}(K''_{+-}(w))\sqcup w^{-1}(K''_{--}(w))=K''\cap K_{-}=K''_{-+}(w)\sqcup K''_{--}(w).$$
 Puisque $w^{-1}$ se restreint en une bijection de $K''$ sur lui-m\^eme, les \'egalit\'es ci-dessus montrent que les nombres d'\'el\'ements de $K''_{+-}(w)$ et de $K''_{-+}(w)$ sont \'egaux. Donc $e''(w)$ est aussi le nombre d'\'el\'ements de $K''_{-+}(w)$. L'application $k\mapsto (\check{\alpha}_{d+1-k,k_{j}^b})_{res}$ est une bijection de $K''_{-+}(w)$ sur l'ensemble des \'el\'ements $\beta$ de $-E$ de la forme $(\check{\alpha}_{k,k_{j}^b})_{res}$ tels que $\beta>0$ et $w^{-1}(\beta)<0$. Finalement, $e'(w)+e''(w)$ est le nombre d'\'el\'ements $\beta\in E\cup(-E)=\check{\Sigma}(i,j)_{res}$ de la forme $(\check{\alpha}_{k,k_{j}^b})_{res}$, tels que $\beta>0$ et $w^{-1}(\beta)<0$. C'est aussi le nombre des $k\in K_{+}$ tels que $w^{-1}(k)\in K_{-}$, c'est-\`a-dire $e_{i}(w)$. La contribution des \'el\'ements de $\check{\Sigma}(i,j)_{res}$ est donc $(-1)^{e_{i}(w)}$.
 
  Le produit sur $i\in I^{-*}$ de ces contributions est $(-1)^{e(w)}$. Il compense exactement  la contribution de $S_{H}(w)S(w)^{-1}$. Le r\'esultat est
 
 pour $j\in I^+\setminus I^{+*}$ et $b=1,2$, ${\bf a}_{j}^b=\chi\circ Norm_{F_{\pm j}/F}$.  
  
 \bigskip
  
  \subsection{Calcul de ${\bf a}_{j}^{\diamond}$ pour $j\in  I^{+*}$}
  
  Soit $j\in I^{+*}$. On doit calculer $a^{\diamond}(w)_{k_{j}}$ pour $w\in W_{F_{j}}$. On a $\rho(w)_{k_{j}}=\chi(w)$.  La restriction de $\chi$ \`a $W_{F_{j}}$ correspond au caract\`ere  $\chi\circ Norm_{F_{j}/F}$ de $F_{j}^{\times}$. On a $(S_{H}(w)S(w)^{-1})_{k_{j}}=(-1)^{e(w)}$. 
  
   Dans le produit 2.9(3), seuls contribuent au $k_{j}$-i\`eme coefficient ceux pour lequel l'indice $j$ est notre $j$. Fixons $i\in I^{-*}$. La contribution des orbites asym\'etriques est un produit de deux termes, dont le premier est un produit de $\check{\beta}(-1)$. On calcule la contribution de ce premier produit comme dans le paragraphe pr\'ec\'edent. Le r\'esultat est le suivant. Notons $\check{\Sigma}(i,j)_{res,asym} $ l'ensemble des \'el\'ements de $ \check{\Sigma}(i,j)_{res}$    dont l'orbite est asym\'etrique. Notons $e_{i,asym}(w)$ le nombre  d'\'el\'ements $\beta\in \check{\Sigma}(i,j)_{res,asym}$ de la forme $(\check{\alpha}_{k,k_{j}})_{res}$, tels que $\beta>0$ et $w^{-1}(\beta)<0$. Alors la contribution est $(-1)^{e_{i,asym}(w)}$.
   
   Soit $\xi\in \Xi(i,j)_{\pm}$, posons $\beta=(\check{\alpha}_{k_{i},\xi(k_{j})})_{res}$, supposons l'orbite de $\beta$ asym\'etrique, calculons la contribution du deuxi\`eme terme de $r_{\beta}(w)^{-1}$.  Pour $n=1,...,N$, la $k_{j}$-i\`eme composante de $\check{\beta}_{n}(\chi_{\beta}(v_{n}(w)))^{-1}$ vaut
 
 (1)  $\chi_{\beta}(v_{n}(w))$ si $\beta_{n}$ est de la forme $(\check{\alpha}_{k,k_{j}})_{res}$, c'est-\`a-dire si $w_{n}^{-1}\xi\in W_{F_{j}}$;
 
 (2) $\chi_{\beta}(v_{n}(w))^{-1}$ si $\beta_{n}$ est de la forme $(\check{\alpha}_{k,d+1-k_{j}})_{res}$, c'est-\`a-dire si $w_{n}^{-1}\xi\in W_{F_{j}}\tau_{j}$;
 
 (3) $1$ sinon.
 
 Notons ${\cal N}$ le sous-ensemble des $n\in \{1,...,N\}$ tels que $w_{n}^{-1}\xi\in W_{F_{j}}$. Pour $n\in {\cal N}$, posons $w'_{n}=\xi^{-1}w_{n}$. L'ensemble $\{w'_{n};n\in {\cal N}\}$ est un ensemble de repr\'esentants du quotient $\xi^{-1}W_{F_{\beta}}\xi\backslash W_{F_{j}}=(\xi^{-1}W_{F_{i}}\xi\cap W_{F_{j}})\backslash W_{F_{j}}$. Pour $n\in {\cal N}$, $v_{n}(w)$ est d\'efini par l'\'egalit\'e $w_{n}w=v_{n}(w)w_{n'}$. On a n\'ecessairement $n'\in {\cal N}$. Posons $v'_{n}(w)=\xi^{-1}v_{n}(w)\xi$. On a $v_{n}(w)\in \xi^{-1}W_{F_{i}}\xi\cap W_{F_{j}}$ et l'\'egalit\'e pr\'ec\'edente \'equivaut \`a $w'_{n}w=v'_{n}(w)w'_{n'}$. La contribution des termes (1) s'\'ecrit
 $$\prod_{n\in {\cal N}}\chi_{\beta}(\xi v'_{n}(w)\xi^{-1}).$$
 C'est la valeur en $w$ du transfert \`a $W_{F_{j}}$ du caract\`ere $w'\mapsto \chi_{\beta}(\xi w'\xi^{-1})$ de $\xi^{-1}W_{F_{i}}\xi\cap W_{F_{j}}$. Comme caract\`ere de $F_{j}^{\times}$, c'est la restriction \`a $F_{j}^{\times}$ du caract\`ere $\chi_{\beta}\circ \xi$ de $\xi^{-1}(F_{\beta}^{\times})$. Pour $\lambda\in F_{j}^{\times}$, on a
 $$\chi_{\beta}\circ \xi(\lambda)=\chi_{i}\circ Norm_{F_{\beta}/F_{i}}\circ \xi(\lambda)=\prod_{\sigma\in \Gamma_{i}/(\Gamma_{i}\cap \xi \Gamma_{j}\xi^{-1})}\chi_{i}(\sigma\xi(\lambda))=\prod_{\sigma\in \Gamma_{i}\xi\Gamma_{j}/\Gamma_{j}}\chi_{i}(\sigma(\lambda)).$$
Un calcul analogue vaut pour la contribution des termes v\'erifiant (2): il suffit d'y remplacer $\xi$ par $\xi\tau_{j}^{-1}$ et les caract\`eres par leurs inverses. On obtient un caract\`ere de $F_{j}^{\times}$ d\'efini par
$$\lambda\mapsto \prod_{\sigma\in \Gamma_{i}\xi\tau_{j}^{-1}\Gamma_{j}/\Gamma_{j}}\chi_{i}(\sigma(\lambda))^{-1}.$$
On se rappelle qu'il est associ\'e \`a $\xi$ un \'el\'ement $\xi'\in \Xi(i,j)$ tel que $\tau_{i}\xi\tau_{j}\in \Gamma_{i}\xi'\Gamma_{j}$. Puisque l'orbite de $\beta$ est asym\'etrique, on a $\xi'\not=\xi$. On a l'\'egalit\'e $\Gamma_{i}\xi\tau_{j}^{-1}\Gamma_{j}=\tau_{i}\Gamma_{i}\xi'\Gamma_{j}$. Puisque $\chi_{i}\circ \tau_{i}=\chi_{i}^{-1}$, le caract\`ere ci-dessus est \'egal \`a 
$$\lambda\mapsto \prod_{\sigma\in \Gamma_{i}\xi'\Gamma_{j}/\Gamma_{j}}\chi_{i}(\sigma(\lambda)).$$
Autrement dit, la contribution des termes (2) est similaire \`a celle des termes v\'erifiant (1): on a seulement remplac\'e $\xi$ par $\xi'$.

Supposons maintenant l'orbite de $\beta$ sym\'etrique. Pour $n=1,...,N$, la $k_{j}$-i\`eme composante de $\check{\beta}_{n}(\chi_{\beta}(v_{n}(w)))^{-1}$ est encore calcul\'ee par (1), (2) et (3) ci-dessus. Notons ${\cal N}_{1}$, resp. ${\cal N}_{2}$, le sous-ensemble des $n\in \{1,...,N\}$ tels que $w_{n}^{-1}\xi\in W_{F_{j}}$, resp. $w_{n}^{-1}\xi\in W_{F_{j}}\tau_{j}$. Remarquons que  $w_{0}\xi\tau_{j}^{-1}\in \xi W_{F_{j}}$. La condition $n\in {\cal N}_{2}$ \'equivaut donc \`a $w_{-n}^{-1}\xi\in W_{F_{j}}$. Pour $n\in {\cal N}_{1}$, posons $w''_{n}=\xi^{-1}w_{n}$. Pour $n\in {\cal N}_{2}$, posons $w''_{n}=\xi^{-1}w_{-n}$. L'ensemble $\{w''_{n}; n\in {\cal N}_{1}\cup {\cal N}_{2}\}$ est un ensemble de repr\'esentants du quotient $(\xi^{-1}W_{F_{i}}\xi\cap W_{F_{j}})\backslash W_{F_{j}}$. Pour $n\in {\cal N}_{1}\cup {\cal N}_{2}$, soient $(n'',v''_{n}(w))\in ({\cal N}_{1}\cup {\cal N}_{2})\times(\xi^{-1}W_{F_{i}}\xi\cap W_{F_{j}})$ tels que $w''_{n}w=v''_{n}(w)w''_{n''}$.  Comparons $v''_{n}(w)$ et $v_{n}(w)$. Rappelons que l'on a une \'egalit\'e $w_{n}w=v_{n}(w)w_{n'}$, avec $n'\in \{\pm 1,...,\pm N\}$. Supposons $n\in {\cal N}_{1}$ et $n'>0$. L'\'el\'ement $n'$ appartient n\'ecessairement \`a ${\cal N}_{1}$ et on a $w''_{n}w=\xi^{-1}v_{n}(w)\xi w''_{n'}$, d'o\`u $v''_{n}(w)=\xi^{-1}v_{n}(w)\xi$. Supposons $n\in {\cal N}_{1}$ et $n'<0$. Alors $-n'$ appartient n\'ecessairement \`a ${\cal N}_{2}$ et on a $w''_{n}w=\xi^{-1}v_{n}(w)\xi w''_{-n'}$, d'o\`u $v''_{n}(w)=\xi^{-1}v_{n}(w)\xi$. Supposons $n\in {\cal N}_{2}$ et $n'>0$. L'\'el\'ement $n'$ appartient n\'ecessairement \`a ${\cal N}_{2}$ et on a $w''_{n}w=\xi^{-1}w_{0}v_{n}(w)w_{0}^{-1}\xi w''_{n'}$, d'o\`u $v''_{n}(w)=\xi^{-1}w_{0}v_{n}(w)w_{0}^{-1}\xi$. Supposons enfin $n\in {\cal N}_{2}$ et $n'<0$. L'\'el\'ement $-n'$ appartient n\'ecessairement \`a ${\cal N}_{1}$ et on a  $w''_{n}w=\xi^{-1}w_{0}v_{n}(w)w_{0}\xi w''_{n'}$, d'o\`u $v''_{n}(w)=\xi^{-1}w_{0}v_{n}(w)w_{0}\xi$. La contribution de $n\in {\cal N}_{1}\cup {\cal N}_{2}$ est donc

$\bullet$ $\chi_{\beta}(\xi v''_{n}(w)\xi^{-1})$ si $n\in {\cal N}_{1}$;

  $\bullet$ $\chi_{\beta}(w_{0}^{-1}\xi v''_{n}(w)\xi^{-1}w_{0})^{-1}$ si $n\in {\cal N}_{2}$ et $n'>0$;
  
  $\bullet$ $\chi_{\beta}(w_{0}^{-1}\xi v''_{n}(w)\xi^{-1}w_{0}^{-1})^{-1}$, si $n\in {\cal N}_{2}$ et $n'<0$.
  
  Les propri\'et\'es de $\chi_{\beta}$ se traduisent en termes galoisiens par $\chi_{\beta}(w_{0}^{-1}vw_{0})=\chi_{\beta}(v)^{-1}$ pour tout $v\in W_{F_{\beta}}$ et $\chi_{\beta}(w_{0}^2)=-1$. Alors le produit des termes ci-dessus est \'egal \`a
  $$(-1)^{e_{\xi}(w)}\prod_{n\in {\cal N}_{1}\cup {\cal N}_{2}}\chi_{\beta}(\xi v''_{n}(w)\xi^{-1}),$$
  o\`u $e_{\xi}(w)$ est le nombre d'\'el\'ements $n\in {\cal N}_{2}$ tels que $n'<0$. Comme dans le cas d'une orbite asym\'etrique, le produit s'interpr\`ete comme un transfert, lequel se traduit par le caract\`ere 
 $$\lambda\mapsto \prod_{\sigma\in \Gamma_{i}\xi\Gamma_{j}/\Gamma_{j}}\chi_{i}(\sigma(\lambda))$$
 de $F_{j}^{\times}$. Par l'application $n\mapsto \beta_{n}$, on voit que $e_{\xi}(w)$ est le nombre de $\beta'$ dans l'orbite de $\beta$, de la forme $(\alpha_{k,k_{j}})_{res}$, tels que $\beta'>0$ et $w^{-1}(\beta')<0$. 
 
 Notons $\check{\Sigma}(i,j)_{res,sym}$ l'ensemble des $\beta\in \check{\Sigma}(i,j)_{res}$ dont la $\Gamma$-orbite est sym\'etrique. Notons $e_{i,sym}(w)$ le nombre d'\'el\'ements $\beta\in \check{\Sigma}(i,j)_{res,sym}$ de la forme $(\check{\alpha}_{k,k_{j}})_{res}$, tels que $\beta>0$ et $w^{-1}(\beta)<0$.

 Le produit sur $\xi\in \Xi(i,j)_{\pm}$ des contributions obtenues est le produit de 
 
 $ \bullet$ $(-1)^{e_{i,sym}(w)}$;  
 
 $\bullet$ la valeur en $w$ du caract\`ere de $W_{F_{j}}$ dont le caract\`ere associ\'e de $F_{j}^{\times}$ est
 $$\lambda\mapsto \prod_{\xi\in \Xi(i,j)}\prod_{\sigma\in \Gamma_{i}\xi\Gamma_{j}/\Gamma_{j}}\chi_{i}(\sigma(\lambda))=\prod_{\sigma\in \Gamma/\Gamma_{j}}\chi_{i}(\sigma(\lambda))=\chi_{i}\circ Norm_{F_{j}/F}(\lambda).$$

 On se rappelle que les orbites asym\'etriques dans $\check{\Sigma}(i,j)_{res}$ avaient d\'ej\`a fourni une contribution $(-1)^{e_{i,asym}(w)}$. La somme $e_{i,sym}(w)+e_{i,asym}(w)$ est le nombre d'\'el\'ements  $\beta\in \check{\Sigma}_{res}(i,j)$ de la forme $(\check{\alpha}_{k,k_{j}})_{res}$, tels que $\beta>0$ et $w^{-1}(\beta)<0$, c'est-\`a-dire $e_{i}(w)$. Quand on fait le produit des signes $(-1)^{e_{i}(w)}$ sur les $i\in I^{-*}$, on obtient $(-1)^{e(w)}$   et ce terme compense le signe qui provient de $S_{H}(w)S(w)^{-1}$. Finalement
 
 pour $j\in I^{+*}$, ${\bf a}^{\diamond}_{j}=\chi\circ Norm_{F_{j}/F}\prod_{i\in I^{-*}}\chi_{i}\circ Norm_{F_{j}/F}$.
 
 \bigskip
 
 \subsection{Fin du calcul de $\Delta_{III}(y,\tilde{x})$}
 
 En utilisant la formule 2.9(1) et nos calculs des diff\'erents caract\`eres qui y interviennent, on obtient que $\Delta_{III}(y,\tilde{x})$ est \'egal au produit des quatre termes suivants:
 
 (1) $\prod_{i\in I^{-*},j\in I^{+*}}\chi_{i}((c_{i}^{-1}\tau_{i}(x_{i}))^{[F_{j}:F]}Norm_{F_{j}/F}(c_{j}\tau_{j}(x_{j})^{-1}))$;
 
 (2) $\prod_{i\in I^{-*}}\chi_{i}(c_{i}^{-1}\tau_{i}(x_{i}))$;
 
 (3) $\prod_{i\in I^{-*}}sgn_{F_{i}/F_{\pm i}}(c_{D}x_{D}^{-1})$;
 
 (4) $\chi(c_{D}x_{D}^{-1}\prod_{j\in I^+}Norm_{F_{j}/F}(c_{j}\tau_{j}(x_{j})^{-1}))$.
 
 \bigskip
  
 \subsection{ Calcul du terme 2.15(3) }
 
  Calculons le d\'eterminant dans $F^{\times}/F^{\times,2}$ de la forme quadratique $\tilde{\theta}$. D'apr\`es la d\'efinition de 1.2, c'est $-\nu$. On doit comparer $\nu$ et le terme $\eta$ d\'efini en 1.6. Ici, le groupe $G_{\theta}$ \'etant sp\'ecial orthogonal impair, tous ses \'epinglages d\'efinis sur $F$ sont conjugu\'es par $G_{\theta}(F)$ et conduisent au m\^eme $\eta$. On peut choisir l'\'epinglage de sorte que l'\'el\'ement nilpotent $N$ associ\'e soit simplement d\'efini par $Ne_{k}=e_{k-1}$ pour $k=2,...,d$ et $Ne_{1}=0$. Alors $\eta=\tilde{\theta}(e_{d},N^{d-1}e_{d})=\tilde{\theta}(e_{d},e_{1})=-\nu$. Consid\'erons la d\'efinition de $c_{D}$ donn\'ee en 2.1. Le d\'eterminant de $\tilde{\theta}$ est le produit de $c_{D}$ et du produit sur tous les $i\in I$ des d\'eterminants des formes quadratiques $(v_{i},v'_{i})\mapsto trace_{F_{i}/F}(\tau_{i}(v_{i})v'_{i}c_{i})$. Fixons $i\in I$. Soit $\delta_{i}\in F_{\pm i}^{\times}$ tel que $F_{i}=F_{\pm i}(\sqrt{\delta_{i}})$ (si $F_{i}$ n'est pas un corps, $\delta_{i}=1$). On v\'erifie que le d\'eterminant de la forme ci-dessus est $Norm_{F_{\pm i}/F}(-\delta_{i})$. L'\'el\'ement $X_{i}$ fix\'e en 2.5 v\'erifie $\tau_{i}(X_{i})=-X_{i}$. Il est donc de la forme $\mu_{i}\sqrt{\delta_{i}}$, avec $\mu_{i}\in F_{\pm i}^{\times}$. Alors $Norm_{F_{i}/F_{\pm i}}(X_{i})=-\mu_{i}^2\delta_{i}$ et le d\'eterminant pr\'ec\'edent est \'egal (modulo un carr\'e) \`a $Norm_{F_{i}/F}(X_{i})$. En utilisant 2.5(2), on obtient
  $$c_{D}\equiv \eta \prod_{i\in I}Norm_{F_{i}/F}((y_{i}-1)(1+y_{i})^{-1})\,\,mod\,\,F^{\times,2}.$$
  Par d\'efinition du polyn\^ome $P_{I}$, on a les \'egalit\'es
  $$P_{I}(-1)=\prod_{i\in I}Norm_{F_{i}/F}(-1-y_{i})=\prod_{i\in I}Norm_{F_{i}/F}(1+y_{i}),$$
  $$P_{I}(1)=\prod_{i\in I}Norm_{F_{i}/F}(1-y_{i})=\prod_{i\in I}Norm_{F_{i}/F}(y_{i}-1).$$
   On en d\'eduit
  
 $$c_{D}\equiv \eta P_{I}(1)P_{I}(-1),\,mod\,\,F^{\times,2}.$$
 Puisque les caract\`eres $sgn_{F_{i}/F_{\pm i}}$ sont d'ordre $2$, le terme 2.15(3) vaut donc
 $$\prod_{i\in I^{-*}}sgn_{F_{i}/F_{\pm i}}(\eta x_{D}P_{I}(1)P_{I}(-1))$$
 
 \bigskip
 
 \subsection{Calcul du terme 2.15(4)}
 
 Soit $j\in I^+$. D'apr\`es 2.5(3), on a 
 $$x_{j}\tau_{j}(x_{j})^{-1}=y_{j}=(1+X_{j})(1-X_{j})^{-1}=(1+X_{j})\tau_{j}(1+X_{j})^{-1} .$$
 Il existe donc $\lambda_{j}\in F_{\pm j}^{\times}$ tel que $x_{j}=\lambda_{j}(1+X_{j})$. Alors
 $$Norm_{F_{j}/F}(c_{j}\tau_{j}(x_{j})^{-1})= Norm_{F_{j}/F}(c_{j}\lambda_{j}^{-1}(1-X_{j})^{-1}).$$
 En utilisant 2.5(2), on obtient
 $$Norm_{F_{j}/F}(c_{j}\tau_{j}(x_{j})^{-1})=Norm_{F_{\pm j}/F}(4^{-1}c_{j}^2\lambda_{j}^{-2})Norm_{F_{j}/F} (1+y_{j})\equiv Norm_{F_{j}/F}(1+y_{j})\,\,mod\,\,F^{\times,2}.$$
 Puisque $\chi$ est d'ordre au plus $2$, le terme 2.15(4) est \'egal \`a
 $$\chi(c_{D}x_{D}\prod_{j\in I^+}Norm_{F_{j}/F}(1+y_{j})).$$
 Comme  en 2.16, le produit intervenant ci-dessus vaut $P_{I^+}(-1)$. En rempla\c{c}ant $c_{D}$ par sa valeur calcul\'ee en 2.16, le terme 2.15(4) vaut
 $$\chi(\eta x_{D}P_{I}(1)P_{I^-}(-1)).$$
 
 \bigskip
 
 \subsection{Disparition des termes 2.8(1) et 2.15(1)}
 
 Fixons $i\in I^{-*}$ et $j\in I^{+*}$. Le produit des contributions de $(i,j)$ \`a 2.8(1) et 2.15(1) est $\chi_{i}(A_{i,j})$ o\`u
 $$A_{i,j}=(c_{i}^{-1}\tau_{i}(x_{i}))^{[F_{j}:F]}Norm_{F_{j}/F}(c_{j}\tau_{j}(x_{j})^{-1})P_{j}(y_{i})Q_{j}(X_{i})^{-1}.$$
 Introduisons la relation d'\'equivalence entre deux \'el\'ements $\mu,\mu'\in \bar{F}^{\times}$: $\mu\equiv_{i}\mu'$ si et seulement si $\mu^{-1}\mu'\in Norm_{F_{i}/F_{\pm i}}(F_{i}^{\times})$. 
 On va montrer que $\chi_{i}(A_{i,j})=1$. Pour cela, il suffit de prouver que $A_{i,j}\equiv_{i} 1$. On a $\tau_{i}(x_{i})^2=\tau_{i}(x_{i})x_{i}y_{i}^{-1}$, d'o\`u
 $$(c_{i}^{-1}\tau_{i}(x_{i}))^{[F_{j}:F]}=(c_{i}^2Norm_{F_{i}/F_{\pm i}}(x_{i}))^{[F_{\pm j}:F]}y_{i}^{-[F_{\pm j}/F]}\equiv_{i} y_{i}^{-[F_{\pm j}:F]}.$$
 On a $Norm_{F_{j}/F}(c_{j})=Norm_{F_{\pm j}/F}(c_{j})^2\equiv_{i} 1$.
  En utilisant la relation 2.5(3) pour $i$ et $j$ on peut \'etablir l'\'egalit\'e suivante:
 $$P_{j}(y_{i})=(1-X_{i})^{-[F_{j}:F]}P_{j}(-1)Q_{j}(X_{i}),$$
 cf. [Li] lemme 7.12.  D'o\`u
 $$A_{i,j}\equiv_{i}  y_{i}^{-[F_{\pm j}:F]}Norm_{F_{j}/F}(\tau_{j}(x_{j}))^{-1}(1-X_{i})^{-[F_{j}:F]}P_{j}(-1).$$
  On a $(1-X_{i})^{-2}=(1-X_{i})^{-1}(1+X_{i})^{-1}y_{i}\equiv_{i} y_{i}$, d'o\`u 
  $$y_{i}^{-[F_{\pm j}:F]}(1-X_{i})^{-[F_{j}:F]}\equiv_{i} 1.$$
  On a aussi 
  $$P_{j}(-1)=Norm_{F_{j}/F}(-1-y_{j})= Norm_{F_{j}/F}((-\tau_{j}(x_{j})-x_{j})\tau_{j}(x_{j})^{-1})$$
  $$=Norm_{F_{\pm j}/F}(x_{j}+\tau_{j}(x_{j}))^2Norm_{F_{j}/F}(\tau_{j}(x_{j}))^{-1}\equiv_{i} Norm_{F_{j}/F}(\tau_{j}(x_{j}))^{-1},$$
  d'o\`u
  $$A_{i,j}\equiv_{i} Norm_{F_{j}/F}(\tau_{j}(x_{j}))^{-2}\equiv_{i} 1.$$
  
  \bigskip
  
  \subsection{Fin du calcul}
  
  Il nous reste \`a calculer le produit de 2.5(1) et des termes 2.8(2) et 2.15(2).  On obtient un produit
  $$\prod_{i\in I^{-*}}sgn_{F_{i}/F_{\pm i}}(\eta c_{i}Q'_{X}(X_{i}))\chi_{i}((y_{i}+1)2^{-1}c_{i}^{-1}\tau_{i}(x_{i})).$$
Soit $i\in I^{-*}$.   Remarquons que $(y_{i}+1)\tau_{i}(x_{i})=(x_{i}+\tau_{i}(x_{i}))$. Ce terme appartient \`a $F_{\pm i}$. Cela permet de remplacer $\chi_{i}$ par $sgn_{F_{i}/F_{\pm i}}$ dans la formule ci-dessus et on obtient
$$(1) \prod_{i\in I^{-*}}sgn_{F_{i}/F_{\pm i}}(B_{i}),$$
o\`u
$$B_{i}=2^{-1}\eta Q'_{X}(X_{i})(y_{i}+1)\tau_{i}(x_{i}).$$
Notons $P_{y}$ le polyn\^ome caract\'eristique de $y$ agissant sur $V^-\oplus V^+$. Les \'egalit\'es 2.5(3) nous permettent d'\'etablir l'\'egalit\'e
  $$2(1-X_{i})^{d-2}P'_{y}(y_{i})=-P_{y}(-1)Q'_{X}(X_{i}),$$
  cf. [Li] lemme 7.12. Les valeurs propres de $y$ agissant dans $V^-\oplus V^+$ sont les $\phi(y_{i})$, pour $i\in I$ et $\phi\in \Phi_{i}$, et $1$ qui intervient une fois (on se rappelle que $(V^-,q^-)$ est un espace quadratique de dimension impaire). D'o\`u $P_{y}(T)=(T-1)P_{I}(T)$, puis $P'_{y}(y_{i})=(y_{i}-1)P'_{I}(y_{i})$ et $-P_{y}(-1)=2P_{I}(-1)$. Comme dans le paragraphe pr\'ec\'edent, $(1-X_{i})^{2}\equiv_{i}y_{i}^{-1}$, d'o\`u 
  $$Q'_{X}(X_{i})\equiv_{i}(1-X_{i})y_{i}^{(3-d)/2}P'_{I}(y_{i}) P_{I}(-1)^{-1}(y_{i}-1).$$
  On a $\tau_{i}(x_{i})=x_{i}\tau_{i}(x_{i})x_{i}^{-1}\equiv_{i}x_{i}^{-1}$. Enfin, gr\^ace \`a 2.5(3),
  $$2(1-X_{i})(y_{i}+1)=4\equiv_{i}1.$$
  On obtient
  $$B_{i}\equiv_{i}\eta x_{i}^{-1} y_{i}^{(3-d)/2}P'_{I}(y_{i})P_{I}(-1)^{-1}(y_{i}-1).$$
  Le facteur de transfert $\Delta_{H,\tilde{G}}(y,\tilde{x})$ est le produit de (1) ci-dessus et des termes calcul\'es en 2.16 et 2.17. On obtient la formule de la proposition 1.10.

 \bigskip
 
 {\bf Bibliographie}

[KS] R. Kottwitz, D. Shelstad: {\it Foundations of twisted endoscopy}, Ast\'erisque 255 (1999)

 [LS] R.P. Langlands, D. Shelstad: {\it On the definition of transfer factors}, Math. Ann. 278 (1987), p.219-271
 
 [Li] W.-W. Li: {\it Transfert d'int\'egrales orbitales pour le groupe m\'etaplectique}, pr\'epublication 2009
 
 [W] J.-L. Waldspurger: {\it Int\'egrales orbitales nilpotentes et endoscopie pour les groupes classiques non ramifi\'es}, Ast\'erisque 269 (2001)

 \bigskip

CNRS-Institut de Math\'ematiques de Jussieu

175, rue du Chevaleret

75013 Paris

e-mail: waldspur@math.jussieu.fr

 \end{document}